\documentclass{amsart}
\usepackage{amssymb, amscd, latexsym}

\newcommand{\bB}{\mathbb{B}}
\newcommand{\bC}{\mathbb{C}}
\newcommand{\bF}{\mathbb{F}}
\newcommand{\bP}{\mathbb{P}}
\newcommand{\bQ}{\mathbb{Q}}
\newcommand{\bR}{\mathbb{R}}
\newcommand{\bZ}{\mathbb{Z}}

\newcommand{\m}{\mathbf{m}}

\newcommand{\cE}{\mathcal{E}}
\newcommand{\cF}{\mathcal{F}}

\newcommand{\cO}{\mathcal{O}}

\newcommand{\ext}{\mathcal{E}xt}

\newcommand{\fC}{\mathfrak{C}}
\newcommand{\fM}{\mathfrak{M}}
\newcommand{\ft}{\mathfrak{t}}

\newcommand{\ch}{\mathrm{ch}}
\newcommand{\diag}{\mathrm{diag}}

\newcommand{\End}{\mathrm{End}}
\newcommand{\Ext}{\mathrm{Ext}_{\cO_X}}
\newcommand{\Ex}{\mathrm{Ext}_{\cO_{X_0}} }

\newcommand{\Hom}{\mathrm{Hom}_{\cO_X}}
\newcommand{\inst}{ {\mathrm{inst}}  }
\newcommand{\pert}{ {\mathrm{pert}}  }
\newcommand{\rank}{\mathrm{rank}}
\newcommand{\td}{\mathrm{td}}
\newcommand{\vir}{\mathrm{vir}}

\newcommand{\tT}{ {\tilde{T}} }

\newcommand{\ep}{\epsilon}

\newcommand{\va}{\vec{a}}

\newcommand{\ve}{ {\vec{e}} }

\newcommand{\vm}{\vec{m}}

\newcommand{\vu}{ {\vec{u}} }
\newcommand{\vD}{ {\vec{D}} }
\newcommand{\vY}{ {\vec{Y}} }

\newcommand{\lz}{ {\ell_0} }
\newcommand{\linf}{ {\ell_\infty} }
\newcommand{\pt}{\mathrm{pt}}
\newcommand{\pa}{\partial}

\newcommand{\MX}{ { \fM_{r,d,n}(X,\linf) } }
\newcommand{\Mx}{ \fM_{r,d,n}(X_0) }
\newcommand{\define}{\stackrel{\mathrm{def}}{=} }
\newcommand{\DY}{(\mathbf{D},\mathbf{Y}) }

\newtheorem{thm}{Theorem}[section]

\newtheorem{fact}[thm]{Fact}
\newtheorem{lm}[thm]{Lemma}
\newtheorem{rem}[thm]{Remark}
\newtheorem{cor}[thm]{Corollary}
\newtheorem{pro}[thm]{Proposition}
\newtheorem{ex}[thm]{Example}
\newtheorem{df}[thm]{Definition}
\newtheorem{nt}[thm]{Notation}

\begin{document}

\title{The Nekrasov Conjecture for Toric Surfaces }
\author{Elizabeth Gasparim}
\address{School of Mathematics, The University of Edinburgh\\
James Clerk Maxwell Building, The King's Buildings, Mayfield Road\\
Edinburgh, EH9 3JZ, Scotland}
\email{Elizabeth.Gasparim@ed.ac.uk}

\author{Chiu-Chu Melissa Liu}
\address{Department of Mathematics,
Columbia University\\
2990 Broadway\\
New York, NY 10027, USA}
\email{ccliu@math.columbia.edu}

\begin{abstract} The
Nekrasov conjecture predicts a relation 
between the partition function for $N=2$ supersymmetric
Yang--Mills theory
and the 
Seiberg-Witten prepotential. 
For instantons on $\mathbb R^4$,
the conjecture  was proved, independently and using
different methods, 
by Nekrasov-Okounkov, Nakajima-Yoshioka, and Braverman-Etingof.
We prove a generalized version of the conjecture for 
instantons on noncompact toric surfaces.
\end{abstract}

\maketitle

\tableofcontents

\section{Introduction}

\subsection{Background}
The Nekrasov conjecture \cite{Ne2} predicts a surprising 
relation between two seemingly unrelated quantities: the  
partition function for $N=2$ supersymmetric Yang--Mills theory,
 defined in terms of instantons on $\bR^4$,  
and the Seiberg-Witten prepotential \cite{SW},  defined 
in terms of period integrals of a family of hyperelliptic
curves.
For gauge group  $U(r)$, 
Nekrasov and Okounkov proved the conjecture for a list of
gauge theories (4d pure gauge theory, 
4d gauge theory with matter, 5d theory compactified on a circle)
\cite{NO},   Nakajima and Yoshioka proved the conjecture for
4d pure gauge theory \cite{NY1} and  
for 5d theory compactified on a circle \cite{NY2} (see also
G\"{o}ttsche-Nakajima-Yoshioka \cite{GNY2}). 
Braverman and Etingof proved the conjecture 
for 4d pure gauge theory with arbitrary gauge groups \cite{BE}.

In this paper we prove a generalized version of the conjecture 
for instantons on noncompact toric surfaces. 
Instantons on toric surfaces have been
studied in \cite{Ne3, GNY1, GNY2}. 

In field theory terms, Nekrasov's insight involves a comparison of the 
infrared and ultraviolet limits of the SUSY gauge theories, as follows.
 The vacuum expectation value of their observables is not sensitive to
 the energy scale. In the ultraviolet, the theory is weakly coupled and 
dominated by instantons; whereas
in the infrared, there appears a relation to the prepotential of the effective 
theory. In this instance, the physical argument is accompanied by completely
rigorous mathematical definitions, thus allowing us to 
prove the conjecture.

\subsection{Partition functions for instantons on noncompact
toric surfaces}
Let  $X_0 = X \setminus \linf$ be an open toric surface that can be  
compactified to a non-singular projective toric surface $X$  by adding a
line at infinity $\linf\cong \bP^1$ with positive self-intersection number,
so that $T_t =(\bC^*)^2$ acts on $X_0$ and on $X$.
Let $\MX$ denote the moduli space of rank $r$ torsion free sheaves  
over $X$ having 
 Chern classes $c_1=d$ and $c_2=n$, and framed over $\linf$.
Then $\MX$ is a 
smooth variety over $\bC$, and it admits a $T_t\times T_e$-action
with isolated fixed points, where $T_e \cong (\bC^*)^r$
is the maximal torus of the complex gauge group $GL(r,\bC)$
which acts on framings.
We define
$$
\int_{\MX}1
$$
by formally applying the Atiyah-Bott localization
formula. The above integral is a rational function in equivariant
parameters $\ep_1,\ep_2\in H^2_{T_t}(\mathrm{pt})$
and $a_1,\ldots,a_r \in H^2_{T_e}(\mathrm{pt})$.
The Nekrasov partition function 
for supersymmetric $SU(r)$ instantons on $X_0$ 
is  defined as 
$$
Z_{X_0,d}^\inst(\ep_1,\ep_2,\va;\Lambda)\define
 \Lambda^{(1-r)d\cdot d}\sum_{n\geq 0} \Lambda^{2rn}
\int_{\MX} 1 
$$
where $\Lambda$ is a formal variable. It lies in 
the ring
$\bQ(\ep_1,\ep_2,a_1,\ldots,a_r)[[\Lambda ]]$.

In further generality, given two multiplicative classes $A,B$ 
we define
$$
Z_{X_0,A,B,d}^\inst(\ep_1,\ep_2,\va;\Lambda)\define
 \Lambda^{(1-r)d\cdot d}\sum_{n\geq 0} \Lambda^{2rn}
\int_{\MX} A_{\tT}(T_{\fM})B_{\tT}(V)
$$ 
where $T_{\fM}$ is the tangent bundle and 
$V$ is the natural bundle on $\MX$ (see Definition \ref{naturalbundle}).

\subsection{Seiberg-Witten prepotential}
We briefly recall the definition of the Seiberg-Witten 
prepotential for 4d pure $SU(r)$ gauge theory. Appendix \ref{sec:SW}
contains a more detailed discussion and definitions for 
other gauge theories. 

Consider the family of hyperelliptic 
curves parametrized by $\Lambda$ and $\vu=(u_2, u_3,\ldots, u_r)$:
$$
C_{\vu} : \Lambda^r\left(w +\frac{1}{w}\right)
=P(z) = z^r + u_2 z^{r-2} + u_3 z^{r-3} + \cdots + u_r.
$$
The parameter space for $\vu$ is called the $\vu$-plane.
The {\em Seiberg-Witten differential}  
$$
dS=\frac{1}{2\pi\sqrt{-1}} z \frac{dw}{w}
$$
is a meromorphic differential  defined on 
the total space of this family such that 
$\displaystyle{ 
\{ \omega_p \define\frac{\partial }{\partial u_p }(dS)\mid  p =2,\ldots,r  \} }$
is a basis of holomorphic differentials on the
genus $(r-1)$ curve $C_\vu$. Choose
a symplectic basis
$\{ A_\alpha, B_\beta\mid \alpha,\beta =2,\ldots,r\}$ 
of $H^1(C_{\vu}, \mathbb Z)$, and define 
$$
a_\alpha = \int_{A_\alpha}dS,\quad  a_{\beta}^D = 2\pi \sqrt{-1} \int_{B_\beta}dS.
$$
Then  the  1-form  $\displaystyle{\sum_{\alpha=2}^r a_\alpha^D d a_\alpha}$ is closed, so
there exists a locally defined function, the {\em Seiberg-Witten
prepotential} $\cF_0$, such that
$$
\sum_{\alpha=2}^r a_\alpha^D d a_\alpha = d \cF_0,\quad
\textup{i.e.,} \quad 
a_{\alpha}^D =  \frac{\partial \cF_0}{\partial a_{\alpha}}.
$$
The above definitions of $dS, a_\alpha, a_\alpha^D$
are the same as those in \cite{NO}, but are $\sqrt{-1}$ times
the corresponding definitions in \cite{NY, NY1}.

\subsection{Nekrasov conjecture}
Let $q_0, q_1$ be the two $T_t$ fixed points
in $\ell_\infty \subset X$, and let
$u,v \in \bZ \ep_1 \oplus \bZ \ep_2$ be
the weights of the $T_t$-action on 
$(N_{\linf/X})_{q_0}$, $(N_{\linf/X})_{q_1}$,
respectively, where $N_{\linf/X}$ is the normal
bundle of $\linf$ in $X$.
If $w$ is the weight of $T_t$-action 
on $T_{q_0}\linf$ and $k =\linf\cdot\linf >0$, then
$$
v = u - kw.
$$
Define
$$
\cF^\inst_{X_0,A,B,d}(\ep_1,\ep_2,\va;\Lambda)
\define - u(u-kw) \log Z^\inst_{X_0,A,B, d}(\ep_1,\ep_2,\va;\Lambda).
$$

We now 
state  the prototype statement of the conjecture for 
toric surfaces,
which will have 8 incarnations.

\bigskip
\paragraph{\bf Main Theorem} (Nekrasov conjecture for toric surfaces:
prototype statement)

{\it \begin{enumerate} 
\item[(a)] $\cF^{\cdots}_{X_0,A,B,d}(\ep_1,\ep_2,\va, {\bf m};\Lambda)$ is analytic
in $\ep_1,\ep_2$ near $\ep_1=\ep_2=0$.
\item[(b)] $\displaystyle{\lim_{\ep_1,\ep_2\to 0} 
\cF^{\cdots}_{X_0,A,B,d}(\ep_1,\ep_2,\va;\Lambda)}
= k\cF_0^{\cdots}(\va,\Lambda)$, 
where
$\cF_0^{\cdots}(\va,\Lambda)$ is the ${\cdots}$part of the 
Seiberg-Witten prepotential  of  matter case $A,B,{\bf m}$,
and $k=\linf\cdot \linf >0$ is the self intersection
number of $\linf$.
\end{enumerate} }

The 8 cases we prove are 
\begin{itemize}
\item {\em   Instanton part: Theorem \ref{thm:Nekrasov-toric}.}
With the  $\phantom{x}^{\cdots}$  replaced by $\phantom{x}^{\inst}$, 
we prove the following cases 
of the conjecture: 
\begin{enumerate}
\item {\em 4d  pure gauge theory}: $A=B=1$, ${\bf m} = \emptyset$.
\item {\em 4d gauge theory with $N_f$ fundamental matter hypermultiplets}:
$A=1$, $B= (E_{\vm})(V)$
is the $T_m$-equivariant Euler class of $V\otimes M$,
where $M$  is the fundamental representation of $U(N_f)$, 
$T_m$ is the maximal torus of $U(N_f)$,
${\bf m} = (m_1, \dots, m_{N_f})$, 
\item {\em 4d gauge theory with one adjoint matter hypermultiplet}:
$A=E_m(T_{\fM})$ is the equivariant Euler class of the 
tangent bundle of the moduli space, $B=1$, ${\bf m} = m$.
\item {\em 5d gauge theory compactified on a circle}: $A= \hat A_\beta (T_{\fM})$
is the $\hat A_\beta$ genus of the tangent bundle (the usual $\hat A$ genus 
being the case $\beta= 1$), $B=1$, $m=\emptyset$ but 
$\cF$ depends on the additional parameter $\beta$.
\end{enumerate}

\item {\em  Perturbative part: Theorem \ref{thm:toric-pert}.}
With the 
$\phantom{x}^{\cdots}$  replaced by $\phantom{x}^{\pert}$, 
we derive 4 more cases 
of the conjecture, with same restrictions as in the first part: 
\begin{enumerate}
\item {\em 4d  pure gauge theory}.
\item {\em 4d gauge theory with $N_f$ fundamental matter hypermultiplets}.
\item {\em 4d gauge theory with one adjoint matter hypermultiplet}.
\item {\em 5d gauge theory compactified on a circle or circumference $\beta$}.
\end{enumerate}
\end{itemize}

The instanton part follows by localization,  from
known results in the $\bC^2$ case. 
 Indeed,
localization calculations yield an expression
of the instanton partition function $Z^\inst_{X_0,A,B,d}$
over $X_0$ in terms of contributions from
vertices ($T_t$ fixed points in $X_0$) and
and from legs ($T_t$ invariant $\bP^1$ in $X_0$).  Each vertex
contributes one copy of the  instanton partition function
of $\bC^2$, for which the singularity along $\ep_1=\ep_2=0$ is
already known.
The contribution from legs does not introduce more poles along
$\ep_1=\ep_2=0$.  A priori, the tangent weights at
all $T_t$ fixed points in $X_0$ appear in the denominator,
but  an argument
similar to that in \cite[Section 6.1]{Ne3}
shows that these poles mostly cancel out, and we are left
with the two normal weights $u$, $u-kw$ at the 
$T_t$ fixed points on $\linf$.
The perturbative part is fairly straightforward.

\subsection{Outline of the paper}
In Section \ref{sec:moduli}, we describe  properties of the instanton moduli spaces.
In Section \ref{sec:torus},  we study  torus actions on these moduli spaces
and the fixed point sets. In Section \ref{sec:Zfunction}, we
introduce a general instanton partition function
depending on two multiplicative classes $A$, $B$ for noncompact toric
surfaces;  different choices of $A$, $B$ give partition functions
of different gauge theories.
Section \ref{sec:instanton} contains localization computations on
instanton moduli spaces, and the proof of the instanton
part of the conjecture. Section \ref{sec:perturbative}
contains definitions of the perturbative
part of the partition function, and the proof
the perturbative part of the conjecture.

\subsection{Acknowledgements} 
This work started during our participation in
 Program for Women and Mathematics at the 
Institute for Advanced Study, Princeton.  We thank the program organizers
Karen Uhlenbeck, Chuu-Lian Terng, Antonella Grassi and Alice Chang
for their encouragement and support. 
We thank Lothar G\"{o}ttsche, Jun Li, Hiraku Nakajima, Nikita Nekrasov,
Tony Pantev, Nathan Seiberg, Constantin Teleman, and Eric Zaslow
for helpful conversations.

\section{Moduli Spaces of Framed Bundles on Surfaces}
\label{sec:moduli}

We work over $\bC$. 
Let $X$ be a non-singular projective surface. Let $\linf \subset X$ 
be a smooth divisor.  In this section, we introduce moduli spaces
of framed bundles on $X$, and describe basic properties of
these moduli spaces,
generalizing the discussion in \cite[Section 2]{NY1} on the
case $X=\bP^2$.
The framed moduli spaces were constructed in much more
setting by Huybrechts-Lehn \cite{HL}.

Given a positive integer $r$, an integer $n$, and
a cohomology class $d \in H^2(X;\bZ)$, let $\MX$ be the moduli
space which parametrizes isomorphism classes of pairs $(E,\Phi)$ such that
\begin{enumerate}
\item $E$ is a torsion free sheaf on $X$ 
which is locally free in a neighborhood of $\linf$.
\item $\rank(E)=r$, $c_1(E)=d$ and $\int_X c_2(E)=n $.
\item $\Phi:  E|_\linf \stackrel{\sim}{\to} \cO_\linf^{\oplus r}$
is an isomorphism called ``framing at infinity''.
\end{enumerate}
Note that (1) and (2) imply $\int_\linf d=0$.  

\subsection{Dimension of the moduli space}

Given a divisor $D\subset X$, 
let $E(-D)=E\otimes \cO_X(-D)$.
\begin{pro}\label{thm:vanish} Suppose that $\linf\cdot \linf >0$. 
\begin{enumerate}
\item[(a)]
For any $(E,\Phi)\in \MX$ we have
$
\Ext^0(E,E(-\linf))=0.
$
\item [(b)] 
Assume in addition that  $\linf \cong \bP^1$. Then
for any $(E,\Phi)\in \MX$ we have
$$
\Ext^0(E,E(-\linf))=\Ext^2(E,E(-\linf))=0.
$$
\end{enumerate}
\end{pro}

\begin{rem}
If $X$ is a non-singular projective surface which 
contains a smooth divisor $\linf\cong \bP^1$ such that 
$k=\linf\cdot \linf >0$.
Then $T_X \bigr|_\linf \cong \cO_{\bP^1}(k)\oplus \cO_{\bP^1}(2)$, so 
$X$ is rationally connected, or equivalently,
$X$ is a rational surface.
The arithmetic genus of $X$ is  $p_a(X)= \chi(\cO_X)-1 =0$. 
\end{rem}

\begin{proof}[Proof of Proposition \ref{thm:vanish}] 
(a) Assuming that $\linf\cdot \linf>0$, we will show that
$$
\Hom(E,E(-\linf))=0.
$$
Let $s$ be a section of $\cO_X(\linf)$ such that
its zero locus is $\linf$. The exact sequence 
$$
0\to E(-(m+1)\linf)\stackrel{s\cdot}{\to} E(-m \linf)\to E(-m\linf)\otimes \cO_D \to 0
$$
induces a long exact sequence
\begin{eqnarray*}
&& 0\to \Hom(E,E(-(m+1)\linf))\to \Hom(E,E(-m\linf)) \\
&& \quad \quad \to \Hom(E,E(-m\linf)\otimes \cO_\linf)\\
&& \quad \to \Ext^1(E,E(-(m+1)\linf)\to \Ext^1(E,E(-m\linf)) \to \cdots
\end{eqnarray*}
where
$$
\Hom(E,E(-m\linf)\otimes \cO_\linf)\cong H^0(\linf, \cO_X(-m\linf)|_\linf)^{\oplus r^2}
$$
since $E|_\linf$ is trivial. Let $k=\linf\cdot \linf>0$. Then
$$
H^0(\linf,\cO_X(-m\linf)|_\linf) \cong H^0(\bP^1, \cO_{\bP^1}(-mk))=0
$$
when $m>0$. So, for any positive integer $m$,
$$
\Hom(E,E(-(m+1)\linf))\to \Hom(E,E(-m\linf))
$$ 
is an isomorphism, and 
$$
\Ext^1(E,E(-(m+1)\linf))\to \Ext^1(E,E(-m\linf))
$$ 
is injective. As a consequence, any
element in $\Hom(E, E(-\linf))$
restricts to zero in a formal neighborhood
of $\linf$ in $X$. So 
$$
\Hom(E,E(-\linf))=0.
$$

\noindent
(b) We now assume that $\linf\cdot \linf>0$ and $\linf\cong \bP^1$.
By  Serre duality,
$\Ext^2(E,E(-\linf))$ is dual to $\Hom(E,E(K_X+\linf))$. 
We will show that
$$
 \Hom(E,E(K_X+\linf))=0.
$$
The exact sequence 
$$
0\to E(K_X-m\linf)\stackrel{s\cdot}{\to} E(K_X+(1-m)\linf)\to E(K_X+(1-m)\linf)\otimes \cO_D \to 0
$$
induces a long exact sequence
\begin{eqnarray*}
&& 0\to \Hom(E,E(K_X-m\linf))\to \Hom(E,E(K_X+(1-m)\linf)) \\
&& \quad\quad \to \Hom(E,E(K_X+(1-m)\linf)\otimes \cO_\linf)\\
&& \quad \to\Ext^1(E,E(K_X-m\linf)\to \Ext^1(E,E(K_X+(1-m)\linf)) \to \cdots \, .
\end{eqnarray*}

$E|_\linf$ is trivial and $K_{\linf}= (K_X + \linf)|_{\linf}$, so
$$
\Hom(E,E(K_X+(1-m)\linf)\otimes \cO_\linf)\cong H^0(\linf, \cO_\linf(K_\linf)
\otimes \cO_X(-m\linf)|_\linf)^{\oplus r^2}.
$$
Note that
$$
H^0(\linf,\cO_\linf(K_\linf)\otimes\cO_X(-m\linf)|_\linf) \cong H^0(\bP^1, \cO_{\bP^1}(-2-mk))=0
$$
for all $m\geq 0$. So, for any nonnegative integer $m$,
$$
\Hom(E,E(K_X-m\linf )\to \Hom(E,E(K_X+(1-m)\linf))
$$ 
is an isomorphism, and 
$$
\Ext^1(E,E(K_X-m\linf))\to \Ext^1(E,E(K_X+(1-m)\linf))
$$ 
is injective. As a consequence, any
element in $\Hom(E, E(K_X+\linf))$
restricts to zero in a formal neighborhood
of $\linf$ in $X$, and we conclude that 
$$
\Hom(E,E(K_X+\linf))=0.
$$
\end{proof}

\begin{cor}
Let $X$ be a non-singular projective surface, and let $\linf$
be a smooth divisor of $X$ such that
$\linf\cdot\linf>0$. Then for any $(E,\Phi)$
in $\MX$,
\begin{eqnarray*}
&& \dim_\bC \Ext^1(E,E(-\linf))
-\dim_\bC \Ext^2(E,E(-\linf)) \\
&&= 2rn + (1-r) d\cdot d -r^2(p_a(X)+p_a(\linf))
\end{eqnarray*}
where $d\cdot d= \int_X d^2$, 
$p_a(X)$ is the arithmetic genus of $X$, and
$p_a(\linf)$ is the arithmetic genus of $\linf$.
\end{cor}

\begin{proof} Let $(E,\Phi)\in \MX$ be locally
free. By Proposition \ref{thm:vanish} (a),
$$
 \dim_\bC \Ext^1(E,E(-\linf))
-\dim_\bC \Ext^2(E,E(-\linf))= -\chi(\End(E)\otimes \cO_X(-\linf)).
$$
Let $\nu \in H^4(X;\bZ)$ be the Poincar\'{e}
dual of $[\mathrm{pt}] \in H_0(X;\bZ)$, and
let $e\in H^2(X;\bZ)$ be the Poincar\'{e}
dual of $[\linf]\in H_2(X;\bZ)$.
By Hirzebruch-Riemann-Roch, 
$$
\chi(\End(E)\otimes \cO_X(-\linf))
=\deg\bigl(\ch(\End(E))\ch(\cO_X(-\linf))\td(T_X)\bigr) .
$$
We have 
\begin{eqnarray*}
&& \ch(\End(E))= \ch(E)\ch(E^\vee)\\
&& = \Bigl(r +d + (\frac{d^2}{2}-n\nu) \Bigr)\Bigl(r -d + (\frac{d^2}{2} -n\nu)\Bigr)
= r^2 + (r-1)d^2-2rn \nu, 
\end{eqnarray*}
$$
\ch(\cO_X(-\linf))= 1 - e +\frac{e^2}{2} = 1-e +\frac{k}{2}\nu, \,\, 
\mbox{for} \,\, k=\linf\cdot \linf >0,
$$
hence
$$
\ch(\End(E))\ch(\cO_X(-\linf))
= r^2 +(-r^2 e) + \Bigl( (r-1)d^2 +( \frac{ kr^2}{2}-2rn)\nu \Big).
$$
We recall that
$$
\td(T_X) =  1+\frac{1}{2} c_1(X)
+\frac{1}{12}(c_1(X)^2 + c_2(X)).
$$
Let $N_{\linf/X}$ be the normal bundle of $\linf$ in $X$. Then
\begin{eqnarray*}
\int_X e c_1(X) &=&\int_{\linf}\bigl( c_1(\linf)+ c_1(N_{\linf/X})\bigr) = 2-2p_a(\linf)+k.
\end{eqnarray*}
Consequently,
\begin{eqnarray*}
&& \deg\left(\ch(\End(E))\ch(\cO_X(-\linf))\td(T_X)\right)\\
&=& \int_X \Bigl(\frac{r^2}{12} (c_1(X)^2+c_2(X))
-\frac{r^2}{2}e c_1(X)
+ (r-1)d^2 +(\frac{kr^2}{2}-2rn)\nu\Bigr)\\
&=& \frac{r^2}{12}\int_X (c_1(X)^2+c_2(X))
-\frac{r^2}{2}(k+2-2p_a(\linf))
+ (r-1)\int_X d^2 + \frac{kr^2}{2}-2rn \\
&=& -2rn + (r-1)\int_X d^2  + r^2 (p_a(X) + p_a(\linf)).
\end{eqnarray*}
\end{proof}

\begin{cor}\label{thm:dim}
Let $X$ be a non-singular projective rational surface, and let $\linf$
be a divisor of $X$ such that $\linf\cong \bP^1$
and $\linf\cdot\linf>0$. Then
$\MX$ is smooth of (complex) dimension
$$
2rn + (1-r) d\cdot d
$$
where $d\cdot d= \int_X d^2$.
\end{cor}

\begin{ex}
Let $X=\bP^2$, and let
$$
\linf=\{ [Z_0,Z_1,Z_2]\in \bP^2 \mid Z_0=0\} \cong \bP^1.
$$
Then $\linf\cdot \linf =1 >0$.
The moduli space
$\fM_{r,d,n}(\bP^2,\linf)$ is nonempty only if
$\int_{\linf} d =0$, which implies $d=0$.
By Corollary \ref{thm:dim}, the moduli space
$\fM_{r,0,n}(\bP^2,\linf)$ is smooth of complex dimension
$2rn$. (See \cite[Corollary 2.2]{NY1}).
\end{ex}

\begin{ex}\label{Fk}
Let $X=\bF_k \define \bP(\cO_{\bP^1}(-k)\oplus \cO_{\bP^1})$ be
the $k^\mathrm{th}$ Hirzebruch surface, where $k$ is a positive
integer. Let 
$$
\lz = \bP(0\oplus \cO_{\bP^1})\cong\bP^1,\quad 
\linf=\bP(\cO_{\bP^2}(-k)\oplus 0)\cong\bP^1.
$$
Then $\lz\cdot \lz= -k <0$ and $\linf\cdot\linf =k>0$.

The moduli space
$\fM_{r,d,n}(\bF_k,\linf)$ is nonempty only if
$\int_{\linf} d =0$, which implies $d = m\lz$
for some $m\in\bZ$. By Corollary \ref{thm:dim}, the moduli space
$\fM_{r,m\lz,n}(\bF_k,\linf)$ is smooth of complex dimension
$2rn + (r-1)k m^2$.
\end{ex}

\begin{ex}
Let $\ell\subset \bP^2$ be a curve of degree 1, and let
$p_1,...,p_k$ be $k$ generic points in $\bP^2$ which are 
disjoint from $\ell$. Let $\pi:\bB_k\to \bP^2$ be the blowup
of $\bP^2$ at $p_1,\ldots,p_k$. Let $\linf=\pi^{-1}(\ell)\cong \bP^1$,
and let $\ell_i =\pi^{-1}(p_i)$ be the exceptional divisors.
Let $e_\infty, e_1,\ldots, e_k \in H^2(\bB_k;\bZ)$
be the Poincar\'{e} duals of 
$[\linf], [\ell_1],\ldots, [\ell_k]$, respectively.
Then 
$$
H^2(\bB_k;\bZ)=\bZ e_\infty \oplus \bZ e_1 \oplus \cdots \bZ e_k.
$$

The moduli space
$\fM_{r,d,n}(\bB_k,\linf)$ is nonempty only if
$\int_{\linf} d =0$, which implies 
$$
d = m_1 e_1 +\cdots + m_k e_k,\quad m_i\in \bZ.
$$
By Corollary \ref{thm:dim}, the moduli space
$\fM_{r,m_1 e_1+\cdots+ m_k e_k ,n}(\bB_k,\linf)$ is smooth of complex dimension
$$
2rn + (r-1)(m_1^2+ \cdots + m_k^2).
$$
\end{ex}

\subsection{The natural bundle}\label{sec:V}

In this subsection, $X$ is a non-singular projective
rational surface, and $\linf$ is a smooth rational curve
in $X$ such that $\linf\cdot \linf >0$.
The proof of the following proposition
is very similar to that of Proposition 
\ref{thm:vanish}.
\begin{pro}\label{thm:vanishV}
$H^0(X, E(-\linf))=H^2(X,E(-\linf)) =0$.
\end{pro}

Let $\cE \to X\times \MX$ be the universal sheaf. Let
$p_1: X\times \MX \to X$ and $p_2: X\times \MX\to \MX$
be the projections to the two factors. 

\begin{df}\label{naturalbundle}
The natural bundle over $\MX$ is  
$$V \define  (R^1p_2)_*(\cE\otimes p_1^*(\cO_X(-\linf))).$$
\end{df}

\begin{cor}\label{thm:rank}
$V$ is a vector bundle of rank 
$$
n-\frac{1}{2} (d\cdot d + c_1(X)\cdot d)
$$ over $\MX$.
\end{cor}
\begin{proof} We use the notation
in the proof of Corollary \ref{thm:dim}.
Let $(E,\Phi)\in \MX$ be locally
free. The rank of $V$ is given by
$
-\chi(E(-\linf)).
$
By Hirzebruch-Riemann-Roch, 
$$
\chi(E(-\linf))
=\deg\left(\ch(E)\ch(\cO_X(-\linf))\td(T_X)\right)
$$
where
\begin{eqnarray*}
\ch(E)&=& r + d +(\frac{d^2}{2}- n\nu)\\
\ch(\cO_X(-\linf))&=& 1 - e +\frac{e^2}{2} = 1-e +\frac{k}{2}\nu
\end{eqnarray*}
$$
\ch(E)\ch(\cO_X(-\linf))
= r +(d-r e) + \Bigl( \frac{d^2}{2} +( \frac{ kr}{2}-n)\nu \Big)
$$
$$
\td(T_X) =  1+\frac{1}{2} c_1(X)
+\frac{1}{12}(c_1(X)^2 + c_2(X)).
$$
Consequently,
\begin{eqnarray*}
&& \deg\left(\ch(E)\ch(\cO_X(-\linf))\td(T_X)\right)\\
&=& \int_X \Bigl(\frac{r}{12} (c_1(X)^2+c_2(X))
+\frac{1}{2}(d-r e) c_1(X)
+ \frac{d^2}{2} +(\frac{kr}{2}-n)\nu\Bigr)\\
&=& \frac{r}{12}\int_X (c_1(X)^2+c_2(X))
-\frac{r}{2}(k+2)
+ \frac{1}{2}\int_X (d^2+ c_1(X)d) + \frac{kr}{2}-n \\
&=& -n + \frac{1}{2}\int_X (d^2+c_1(X)d)  + r p_a(X)
\end{eqnarray*}
where $p_a(X)=0$ since $X$ is a rational surface.
\end{proof}

\section{Torus Action and Fixed Points}
\label{sec:torus}

In this section, $X$ is a non-singular projective toric surface.
Therefore $T_t\define (\bC^*)^2$ acts on $X$. 
We use notation similar to that in \cite[Section 2, 3]{NY1}.

\subsection{Torus action on the surface}
We assume that $\linf$ is a $T_t$-invariant $\bP^1$ in $X$, and
$\linf \cdot \linf = k>0$.   Then $X_0=X\setminus \linf$ is
a non-singular, quasi-projective toric surface.  Let $\Gamma$ be a graph
such that the vertices of $\Gamma$ are in one-to-one
correspondence with the $T_t$ fixed points in $X_0$,
and two vertices are connected by an edge if and
only if the corresponding fixed points are connected
by a $T_t$-invariant $\bP^1$. Then 
$\Gamma$ is a chain, so
$\# V(\Gamma)- \# E(\Gamma)=1$, and
$$
\chi(X_0)= \# V(\Gamma) =\chi(X)-2,
$$
where $E(\Gamma)$ is the set of edges in $\Gamma$
and $V(\Gamma)$ is the set of vertices in $\Gamma$.
Let $p_v$ be the $T_t$ fixed point in $X_0$
which corresponds to $v\in V(\Gamma)$, and
let $\ell_e$ be the $T_t$-invariant $\bP^1$ which
corresponds to $e\in E(\Gamma)$. Any $T_t$-invariant
divisor $D$ in $X$ disjoint from $\linf$ is of the form
$$
D= \sum_{e\in E(\Gamma)} m_e \ell_e \cong H_2(X_0;\bZ)
$$
where $m_e\in\bZ$.

\subsection{Torus action on moduli spaces}
Let $T_e$ be the maximal torus of $GL(r,\bC)$ consisting of diagonal
matrices, and let $\tT = T_t \times T_e$. We define an action
of $\tT$ on $\MX$ as follows: for 
$(t_1,t_2)\in T_t$, let $F_{t_1,t_2}$ be the automorphism of 
$X$ defined by 
$F_{t_1,t_2}(x) = (t_1,t_2)\cdot x$.
Given $\ve=\diag(e_1,\ldots,e_r)\in T_e$, let $G_\ve$
denote the isomorphism of $\cO_\linf^{\oplus r}$ given by
$(s_1,\ldots,s_r)\mapsto (e_1 s_1,\ldots, e_r s_r)$.
For $(E,\Phi)\in \MX$, we define
$$
    (t_1,t_2,\ve)\cdot (E,\Phi)
    = \left((F_{t_1,t_2}^{-1})^* E, \Phi'\right),
$$
where $\Phi'$ is the composite of homomorphisms
$$
   (F_{t_1,t_2}^{-1})^* E|_\linf 
   \xrightarrow{(F_{t_1,t_2}^{-1})^*\Phi}
   (F_{t_1,t_2}^{-1})^* \cO_\linf^{\oplus r}
   \stackrel{\phi_{t_1,t_2} }{\longrightarrow} \cO_\linf^{\oplus r}
   \xrightarrow{G_\ve} \cO_\linf^{\oplus r}.
$$
Here $\phi_{t_1,t_2}$ is the homomorphism given by the action.

\subsection{Torus fixed points in moduli spaces}
The fixed points set $\MX^{\tT}$ consists of
\( 
   (E,\Phi) = (I_1(D_1),\Phi_1) \oplus \cdots \oplus(I_2(D_r),\Phi_r)
\)
such that
\begin{enumerate}
\item $I_\alpha(D_\alpha)$ is a tensor product 
$I_\alpha\otimes\cO_X(D_\alpha)$, where
$D_\alpha$ is a $T_t$-invariant divisor which does not 
intersect $\linf$, and $I_\alpha$ is  the ideal sheaf of a $0$-dimensional 
subschemes $Q_\alpha$ contained in $X_0$.

\item $\Phi_\alpha$ is an isomorphism from $(I_\alpha)_\linf$ 
to the $\alpha$th factor of $\cO_{\ell_\infty}^{\oplus r}$.
\item $I_\alpha$ is fixed by the action of $T_t$.
\end{enumerate}

The support of $Q_\alpha$ must be contained in
$X_0^{T_t}$, the $T_t$ fixed points set of $X_0$.
Thus $Q_\alpha$ is a union of $\{ Q_\alpha^v \mid v\in V(\Gamma)\}$ 
where $Q_\alpha^v$ is a subscheme supported at the $T_t$-
fixed point $p_v\in X_0$. If we take a coordinate system
$(x,y)$ around $p_v$, then the ideal of $Q_\alpha^v$
is generated by monomials $x^i y^j$,
So $Q_\alpha^v$ corresponds to a Young diagram $Y_\alpha^v$.

Therefore the fixed point set is parametrized by $2r$-tuples
$$
\DY=(D_1,\vY_1,\ldots, D_r,\vY_r)
$$
where
$$
D_\alpha \in \bigoplus_{e\in E(\Gamma)} \bZ \ell_e \cong H_2(X_0;\bZ),\quad
\vY_\alpha = \{ Y_\alpha^v \mid v\in V(\Gamma) \},
$$
and each $Y_\alpha^v$ is a Young diagram.
Let
$$
|\vY_\alpha|=\sum_{v\in V(\Gamma)} |Y_\alpha^v|.
$$

Let $d^\vee \in H_2(X;\bZ)$ be the Poincar\'{e} dual
of $d\in H^2(X;\bZ)$. Then $\int_\linf d=0$ implies
$$
d^\vee \in \bigoplus_{e\in E(\Gamma)} \bZ [\ell_e].
$$

 The constraints are
\begin{equation}\label{eqn:cone}
\sum_{\alpha}D_\alpha= d^\vee,
\end{equation}
\begin{equation} \label{eqn:ctwo}
\sum_{\alpha=1}^r |\vY_\alpha| + \sum_{\alpha<\beta} D_\alpha\cdot  D_\beta =n .
\end{equation}
Note that
$$
2r \sum_{\alpha<\beta}D_\alpha \cdot D_\beta+(1-r) d^\vee \cdot d^\vee
=  (1-r)\sum_\alpha D_\alpha^2  + 2 \sum_{\alpha<\beta} D_\alpha\cdot D_\beta
= -\sum_{\alpha<\beta} (D_\alpha-D_\beta)^2,
$$
so \eqref{eqn:ctwo} can be rewritten as
\begin{equation}
2r\sum_{\alpha=1}^r |\vY_\alpha|
-\sum_{\alpha<\beta}(D_\alpha-D_\beta)^2
=2rn + (1-r) d\cdot d = \dim_\bC \MX. 
\end{equation}

\section{Gauge Theory Partition Functions}
\label{sec:Zfunction}

We refer to Appendix \ref{sec:equivariant}
for a brief review of equivariant cohomology
and integration of an equivariant cohomology class over
a possibly non-compact manifold. 

\subsection{Equivariant parameters}
For $i=1,2$, let
$p_i\colon BT_t \cong \bP^\infty \times \bP^\infty$
be the projection to the $i$-th factor, 
and let $\ep_i = (c_1)_{T_t}(p_i^*\cO(1))$. Then 
$$
H_{T_t}^*(\pt;\bQ)= H^*(BT_t;\bQ)=\bQ[\ep_1,\ep_2].
$$
Let $t_i = e^{\ep_i} = \ch_1(p_i^*\cO(1))$.

Similarly, for $j=1,\ldots,r$, let 
$q_j: BT_e \cong (\bP^\infty)^r \to \bP^\infty$
be the projection to the $j$-th factor, 
and let $a_j = (c_1)_{T_t}(q_j^*\cO(1))$. Then
$$
H^*_{T_e}(\pt;\bQ)= H^*(BT_e;\bQ) =\bQ[a_1,\ldots,a_r].
$$
Let $e_j = e^{a_j}=\ch_1(q_j^*\cO(1))$.
We write $\va=(a_1,\ldots,a_r)$ and
$\ve=(e_1,\ldots,e_r)=(e^{a_1},\ldots, e^{a_r})$.

\subsection{Multiplicative classes of the tangent and natural bundles}

Recall that a multiplicative class $c$ is a characteristic class which
satisfies $c(E_1 \oplus E_2) = c(E_1)c(E_2)$. Such a class
is determined by a formal power series $f(x)$ satisfying
$c(L) = f(c_1(L))$ for a line bundle $L$ and
$c(E)=f(x_1) \cdots f(x_r)$ where $x_1,\ldots,x_r$ are Chern
roots of $E$.

Let $A$, $B$ be multiplicative classes associated
to formal power series $f(x)$, $g(x)$, respectively.  Then
$$
\int_{\MX} A_{\tT}(T_{\fM})B_{\tT}(V)
\in \bQ[[\ep_1,\ep_2,\va ]]_\m \subset \bQ(( \ep_1,\ep_2,\va)),
$$ 
where $T_{\fM}$ is the tangent bundle of $\MX$, $V$ is defined
in Definition \ref{naturalbundle}, and 
$\bQ[[\ep_1,\ep_2,\va]]_\m$ is the localization of 
the ring $\bQ[[\ep_1,\ep_2,\va]]$ at the maximal ideal $\m$ generated
by $\ep_1,\ep_2,a_1,\ldots,a_r$. If
$f(x)$ and $g(x)$ are polynomials, then
$$
\int_{\MX} A_{\tT}(T_{\fM})B_{\tT}(V)
\in \bQ[\ep_1,\ep_2,\va ]_\m
\subset \bQ(\ep_1,\ep_2,\va).
$$

Let $X_0=X\setminus \linf$. 
Given 
$$
d\in \{ \gamma \in H^2(X;\bZ)\mid \int_\linf \gamma =0\} \cong H^2_c(X_0;\bZ)
$$
let $d^\vee \in H_2(X;\bZ)$ be its Poincar\'{e} dual.
(Here $H_c^*$ is the compact cohomology.) Then
$$
d^\vee \in \bigoplus_{e\in E(\Gamma)}\bZ \ell_e \cong H_2(X_0;\bZ).
$$
We define
\begin{eqnarray*}
&& Z_{X_0,A,B,d}^{\inst}(\ep_1,\ep_2,\va;\Lambda) \\
&=& \sum_{n\geq 0}\Lambda^{\dim_\bC \MX}
\int_{\MX}A_{\tT}(T_{\fM})B_{\tT}(V) \\
&=& \Lambda^{(1-r)d\cdot d}\sum_{n\geq 0} \Lambda^{2rn}
\int_{\MX}A_{\tT}(T_{\fM}) B_{\tT}(V)\\
&=&\sum_{\sum D_\alpha= d^\vee} 
\Lambda^{-\sum_{\alpha<\beta} (D_\alpha-D_\beta)^2}
\sum_{\vY_\alpha} \Lambda^{\sum_\alpha |\vY_\alpha|}
\frac{A_{\tT}(T_{\DY}\MX) B_{\tT}(V_{\DY})}{e_{\tT}(T_{\DY}\MX)}\\
&=&\sum_{\sum D_\alpha= d^\vee} \sum_{\vY_\alpha}
\prod(\Lambda\frac{f(x_i)}{x_i})\prod g(y_j) \in \bQ((\ep_1,\ep_2,\va))[[\Lambda]]
\end{eqnarray*}
where $d^\vee$ $x_i$ are $\tT$-equivariant Chern roots of $T_{\DY}\MX$ and
$y_j$ are $\tT$-equivariant Chern roots of $V_{\DY}$. If $f(x)$, $g(x)$
are polynomials then
$$
Z_{X_0,A,B,d}^{\inst}(\ep_1, \ep_2,\va;\Lambda)\in\bQ(\ep_1,\ep_2,\va)[[\Lambda]].
$$

Sometimes we allow $A$ and $B$ to
depend on extra parameters, then 
$Z_{X,A,B,d}^\inst$ will depend on extra parameters as well.

Introduce variables $\{ Q_e \mid e\in E(\Gamma)\}$. 
Given $d\in H^2_c (X_0;\bZ)$, define
$$
Q^d =\prod_{e\in E(\Gamma)} Q_e^{\int_{\ell_e} d}.
$$
We define a generating function
\begin{eqnarray*}
&& Z_{X_0, A, B}^{\inst}(\ep_1,\ep_2,\va;\Lambda,Q)
\define \sum_{d\in H^2_c(X_0;\bZ)}
Q^d Z_{X_0,A,B,d}^{\inst}(\ep_1,\ep_2,\va;\Lambda)\\
&&= \sum_{d\in H^2_c(X_0;\bZ)}\sum_{n\geq 0}
Q^d \Lambda^{(1-r)d\cdot d + 2rn}\int_{\MX} 
A_{\tT}(T_{\fM})B_\tT(V).
\end{eqnarray*}

\subsection{4d pure gauge theory}
\label{sec:four-pure}
Nekrasov instanton partition functions of 4d pure gauge theory are given by 
\begin{eqnarray*}
Z_{X_0,d}^\inst(\ep_1,\ep_2,\va;\Lambda)\define
 \Lambda^{(1-r)d\cdot d}\sum_{n\geq 0} \Lambda^{2rn}
\int_{\MX} 1,\\
Z_{X_0}^\inst(\ep_1,\ep_2,\va;\Lambda,Q)\define
\sum_{d\in H^2_c(X_0;\bZ) } Q^d Z_{X_0,d}^\inst(\ep_1,\ep_2,\va;\Lambda).
\end{eqnarray*}
We have
\begin{eqnarray*}
Z_{X_0,d}^\inst(\ep_1,\ep_2,\va;\Lambda)&=&
Z_{X_0,A=1,B=1,d}^{\inst}(\ep_1,\ep_2,\va ;\Lambda) ,\\
Z_{X_0}^\inst(\ep_1,\ep_2,\va;\Lambda,Q)&=&
Z_{X_0,A=1, B=1}^{\inst}(\ep_1,\ep_2,\va;\Lambda,Q).
\end{eqnarray*}

We define a grading on the ring $\bQ((\ep_1,\ep_2,\va))[[\Lambda]]$ by
$$
\deg \Lambda  =\deg \ep_1 =\deg \ep_2 =\deg a_\alpha =2.
$$
Then $Z_{X_0,d}^\inst(\ep_1,\ep_2,\va;\Lambda)
\in \bQ((\ep_1,\ep_2,\va))[[\Lambda]]$ is homogeneous of degree 0.

\subsection{4d gauge theory with $N_f$ fundamental matter hypermultiplets}
\label{sec:four-fundamental}
Let $T_m$ be the maximal torus of $U(N_f)$. Then $H^*_{T_m}(\pt)\cong \bQ[m_1,\ldots, m_{N_f}]$.
Let $M$ be the fundamental representation of $U(N_f)$, and write $\vm=(m_1,\ldots,m_{N_f})$. 
Let $V$ be the natural vector bundle as in Definition \ref{naturalbundle}; it is a $\tT$-equivariant
vector bundle over $\MX$.

Nekrasov instanton partition functions of 4d gauge theory with $N_f$ fundamental matter hypermultiplets
are given by 
\begin{eqnarray*}
Z_{X_0,d}^\inst(\ep_1,\ep_2,\va,\vm;\Lambda)
&\define& \Lambda^{(1-r)d\cdot d}\sum_{n\geq 0} \Lambda^{2rn}
\int_{\MX} (c_{top})_{\tT\times T_m}(V\otimes M)\\
&=& \Lambda^{(1-r)d\cdot d}\sum_{n\geq 0} \Lambda^{2rn}
\int_{\MX} \prod_{f=1}^{N_f}(E_{m_f})_{\tT}(V)
\end{eqnarray*}
where $E_t$ is the multiplicative class associated to $f(x)= t+x$, so that
$$
E_t(V)= t^k + c_1(V) t^{k-1}+\cdots+ c_n(V),\quad k= \rank_\bC V.
$$
$$
Z_{X_0}^\inst(\ep_1,\ep_2,\va,\vm;\Lambda,Q)\define
\sum_{d\in H^2_c(X_0;\bZ)} Q^d Z_{X_0,d}^\inst(\ep_1,\ep_2,\va,\vm;\Lambda).
$$

Let $E_{\vm}=\prod_{f=1}^{N_f}E_{m_f}$. Then
\begin{eqnarray*}
Z_{X_0,d}^\inst(\ep_1,\ep_2,\va,\vm;\Lambda)
&=&Z^\inst_{X_0,A=1,B=E_{\vm},d }(\ep_1,\ep_2,\va;\Lambda)\\
Z_{X_0}^\inst(\ep_1,\ep_2,\va,\vm;\Lambda,Q)
&=& Z^\inst_{X_0, A=1,B=E_{\vm} }(\ep_1,\ep_2,\va;\Lambda,Q).
\end{eqnarray*}

\subsection{4d gauge theory with one adjoint matter hypermultiplet}
\label{sec:four-adjoint}
Nekrasov  instanton partition functions of 4d gauge theory with one
adjoint matter hypermultiplet
are given by 
\begin{eqnarray*}
Z_{X_0,d}^\inst(\ep_1,\ep_2,\va,m;\Lambda)
&\define&  \Lambda^{(1-r)d\cdot d}\sum_{n\geq 0} \Lambda^{2rn}
\int_{\MX} (E_m)_{\tT}(T_\fM)\\
Z_{X_0}^\inst(\ep_1,\ep_2,\va,m;\Lambda,Q)&\define&
\sum_{d\in H^2_c(X_0;\bZ)} Q^d Z_{X_0,d}^\inst(\ep_1,\ep_2,\va,m;\Lambda).
\end{eqnarray*}
We have
\begin{eqnarray*}
Z_{X_0,d}^\inst(\ep_1,\ep_2,\va,m;\Lambda)
&=&Z^\inst_{X_0, A=E_m,B=1, d }(\ep_1,\ep_2,\va;\Lambda)\\
Z_{X_0}^\inst(\ep_1,\ep_2,\va,m;\Lambda,Q)
&=&Z^\inst_{X_0, A=E_m,B=1}(\ep_1,\ep_2,\va;\Lambda,Q).
\end{eqnarray*}

\subsection{5d gauge theory compactified on a circle of circumference $\beta$}
\label{sec:five}

Let $\widehat{A}_\beta$ be the
multiplicative class associated to $f_\beta(x)=\displaystyle{ \frac{\beta x/2}{\sinh(\beta x/2)} }$.
For a complex vector bundle $E$,  $\widehat{A}_1 (E)=\widehat{A}(E)$ is the $\widehat{A}$-genus
of $E$. The index of the Dirac operator on a complex manifold $M$ is given by
$$
\int_M \widehat{A}(T_M).
$$

The Nekrasov partition functions of 5d gauge theory compactified on
a circle of circumference $\beta$ are given by
\begin{eqnarray*}
&& Z^{\inst}_{X_0,d}(\ep_1,\ep_2,\va;\Lambda,\beta)=
\Lambda^{(1-r)d\cdot d}\sum_{n\geq 0} \Lambda^{2rn}
\int_{\MX}(\widehat{A}_\beta)_{\tT}(T_\fM),\\
&& Z^{\inst,(m)}_{X_0}(\ep_1,\ep_2,\va;\Lambda,Q,\beta)=
\sum_{d\in H^2_c(X_0;\bZ)} Q^d  Z^{\inst}_{X_0,d}(\ep_1,\ep_2,\va;\Lambda,\beta)
.\end{eqnarray*}
We have
\begin{eqnarray*}
Z_{X_0,d}^{\inst}(\ep_1,\ep_2,\va;\Lambda,\beta) 
&=& Z^{\inst}_{X_0,A=\widehat{A}_\beta, B=1,d}(\ep_1,\ep_2,\va;\Lambda),\\
Z_{X_0}^{\inst}(\ep_1,\ep_2,\va;\Lambda,Q,\beta) 
&=& Z^{\inst}_{X_0, A=\widehat{A}_\beta,B=1}(\ep_1,\ep_2,\va;\Lambda,Q).
\end{eqnarray*}
Note that $\displaystyle{\lim_{\beta\to 0} f_\beta(x)=1}$, so the 
partition function of 5d gauge theory compactified on a circle of 
circumference $\beta$ specializes to the one
of 4d pure gauge theory
as $\beta\to 0$, that is:
\begin{eqnarray*}
\lim_{\beta\to 0} Z_{X_0,d}^{\inst}(\ep_1,\ep_2,\va;\Lambda,\beta)
&=&Z_{X_0,d}^{\inst}(\ep_1,\ep_2,\va;\Lambda),\\
\lim_{\beta\to 0} Z_{X_0}^{\inst}(\ep_1,\ep_2,\va;\Lambda,Q,\beta)
&=&Z_{X_0}^{\inst}(\ep_1,\ep_2,\va;\Lambda,Q).
\end{eqnarray*}

\subsection{Hirzebruch $\chi_y$ genus}\label{sec:Hirzebruch}
Let 
$$
(\chi_y)_\tT(\MX)=\sum_{p=0}^N (-y)^p
\sum_{q=0}^N (-1)^q \ch_{\tT} H^q(\MX, \Lambda^p T^* \MX)
$$
be the $\tT$-equivariant Hirzebruch $\chi_y$ genus,
where $N=\dim_\bC \MX$.  In particular,
$$
(\chi_0)_\tT(\MX)=\chi_\tT(\MX, \cO).
$$
By Hirzebruch-Riemann-Roch,
$$
(\chi_y)_\tT(\MX)=\sum_{p=0}^N (-y)^p
\int_{\MX}\td_\tT(\fM)\ch_\tT(\Lambda^p T^*\fM)
$$
where $\fM=\MX$. Define
\begin{eqnarray*}
Z^\inst_{X_0,d}(\ep_1,\ep_2,\va;\Lambda,y)
&=&\Lambda^{(1-r)d\cdot d}\sum_{n\geq 0} \Lambda^{2rn}
(\chi_y)_\tT(\MX),\\
Z^\inst_{X_0}(\ep_1,\ep_2,\va;\Lambda,Q,y)
&=&\sum_{d\in H^2_c(X_0;\bZ)}Q^d Z^\inst_{X_0,d}(\ep_1, \ep_2,\va;\Lambda, y)
.\end{eqnarray*}
Then
\begin{eqnarray*}
Z^\inst_{X_0,d}(\ep_1,\ep_2,\va;\Lambda,y)
&=&Z^\inst_{X_0,A=A_y,B=1,d}(\ep_1,\ep_2,\va;\Lambda)\\
Z^\inst_{X_0}(\ep_1,\ep_2,\va;\Lambda,Q,y)
&=&Z^\inst_{X_0,A=A_y,B=1}(\ep_1,\ep_2,\va;\Lambda,Q),
\end{eqnarray*}
where $A_y$ is the multiplicative class
associated to 
$$
f_y(x)= \frac{x(1-y e^{-x})}{1-e^{-x}}.
$$
In particular, 
$$
f_0(x)=\frac{x}{1-e^{-x}},\quad
f_1(x)=x,
$$
so $A_0(E)=\td(E)$ and $A_1(E)=e(E)$.

\subsection{Elliptic genus}\label{sec:elliptic}
Let $A_{y,q}$ be the  multiplicative
class associated to 
$$
y^{-1/2} x\prod_{n\geq 1}\frac{(1-y q^{n-1} e^{-x})(1-y^{-1} q^n e^x)}
{(1-q^{n-1} e^{-x})(1-q^n e^x)}
$$
The $\tT$-equivariant elliptic genus of $\fM$ is given by
$$
\chi_\tT (\MX, y,q)
=\int_{\MX} A_{y,q}(T_{\fM}).
$$
Define
\begin{eqnarray*}
& & Z_{X_0,d}^\inst(\ep_1,\ep_2,\va;\Lambda,y,q)
\define  Z_{X_0,A=A_{y,q},B=1,d}^\inst(\ep_1,\ep_2,\va;\Lambda)\\
&&=\Lambda^{(1-r)d\cdot d} \sum_{n\geq 0} \Lambda^{2rn}\chi_\tT(\MX,y,q),\\
& & Z_{X_0}^\inst(\ep_1,\ep_2,\va;\Lambda,Q,y,q)
\define  Z_{X_0,A=A_{y,q},B=1}^\inst(\ep_1,\ep_2,\va;\Lambda,Q)\\
&&=\sum_{d\in H^2_c(X_0,\bZ)}\sum_{n\geq 0}
Q^d \Lambda^{(1-r)d\cdot d+ 2rn}\chi_\tT(\MX,y,q).
\end{eqnarray*}

\section{The Instanton Part}
\label{sec:instanton}

In this section, we calculate the partition functions defined
in Section \ref{sec:Zfunction}.

\subsection{The tangent bundle: adjoint representation} \label{sec:Tinstanton}
Let $(E,\Phi)\in \MX$ be a fixed point of $\tT$-action corresponding to
$\DY=(D_1,\vY_1,\ldots,D_r, \vY_r)$. 
We want to compute
$$
   \ch_{\tT}T_{(E,\Phi)} \MX = \ch_{\tT} \Ext^1(E,E(-\linf))
=-\ch_{\tT} \Ext^*(E,E(-\linf)).
$$ 
We have
$$
E= I_1(D_1)\oplus \cdots \oplus I_r(D_r),
$$
so 
\begin{eqnarray*}
-\ch_{\tT}\Ext^*(E,E(-\linf)) 
&=& -\sum_{\alpha,\beta}\ch_{\tT}\Ext^*(I_\alpha(D_\alpha), I_\beta(D_\beta-\linf))\\
&=& -\sum_{\alpha,\beta}e^{a_\beta-a_\alpha}
\ch_{T_t}\Ext^*(I_\alpha(D_\alpha), I_\beta(D_\beta-\linf)).
\end{eqnarray*}
Let
\begin{eqnarray*}
L_{\alpha,\beta}(t_1,t_2)&=& -\ch_{T_t}\Ext^*(\cO_X(D_\alpha), \cO_X(D_\beta-\linf))\\
&=& -\chi_{T_t}(X, \cO_X(D_\beta-D_\alpha-\linf))
\end{eqnarray*}
\begin{eqnarray*}
M_{\alpha,\beta}(t_1,t_2)&=& \ch_{T_t}\Ext^*(\cO_X(D_\alpha),\cO_X(D_\beta-\linf))\\
&& -\ch_{T_t}\Ext^*(I_\alpha(D_\alpha),I_\beta(D_\beta-\linf)) .
\end{eqnarray*}
Then
\begin{equation}\label{eqn:T-ML}
\ch_\tT T_{(E,\Phi)}\MX =\sum_{\alpha,\beta=1}^r e^{a_\beta-a_\alpha}
\left( M_{\alpha,\beta}(t_1,t_2) + L_{\alpha,\beta}(t_1,t_2)\right).
\end{equation}
So it remains to compute $M_{\alpha,\beta}(t_1,t_2)$ and $L_{\alpha,\beta}(t_1,t_2)$.

\subsubsection{$M_{\alpha,\beta}(t_1,t_2)$}\label{sec:M-T}
Let $\chi_{D_\alpha}^v\in \Hom(T_t,\bC^*)$ be the characters  of the
$T_t$-equivariant line bundle $\cO_X(D_\alpha)$
at the $T_t$ fixed point $p_v\in X_0$, 
and let $\chi_1^v, \chi_2^v\in \Hom(T_t,\bC^*)$ be the characters of $T_{p_v}X$. 
Then $\chi_{D_\alpha}^v, \chi_1^v, \chi_2^v$ are monomials in
$t_1, t_2$.

Let $\ft_t$ be the Lie algebra of $T_t$.
Define weights $w_{D_\alpha}^v, w_1^v, w_2^v \in \Hom(\ft_t,\bC)= \ft_t^\vee$ by
$$
e^{w_{D_\alpha}^v } =\chi_{D_\alpha}^v,\quad
e^{w_1^v}=\chi_1^v,\quad e^{w_2^v}=\chi_2^v.
$$

Given a partition (Young diagram) $S$ and a box $s\in S$, 
let $a_S(s)$ and $l_S(s)$ be the arm-length
and leg-length of $s$ (see e.g. \cite[Figure 2]{NY}).
Given two partitions $S$, $T$, let
\begin{equation}\label{eqn:MST}
M_{S,T}(t_1,t_2)=\sum_{s\in S} t_1^{-l_T(s)} t_2^{a_S(s)+1}
+\sum_{t\in T} t_1^{l_S(t)+1} t_2^{-a_T(t)}
\end{equation}
\begin{equation}
\begin{aligned}
N_{S,T}(\ep_1,\ep_2)& \define M_{S,T}(e^{\ep_1}, e^{\ep_2}) \\
&=\sum_{s\in S} e^{-l_T(s) \ep_1 +(a_S(s)+1)\ep_2}
+\sum_{t\in T} e^{(l_S(t)+1)\ep_1 -a_T(t)\ep_2}.
\end{aligned}
\end{equation}
The expression \eqref{eqn:MST} was introduced in  \cite[Equation (4.45)]{FP}.
(See also \cite[Lemma 3.2]{EG} and \cite[Theorem 2.1]{NY1}.)

\begin{pro}[vertex contribution to the tangent bundle]\label{thm:T-M} 
\begin{eqnarray*}
M_{\alpha,\beta}(t_1,t_2)&=&\sum_{v\in V(\Gamma)}
\frac{\chi_{D_\beta}^v(t_1,t_2)}{\chi_{D_\alpha}^v(t_1,t_2)}
M_{Y_\alpha^v, Y_\beta^v}(\chi_1^v(t_1,t_2), \chi_2^v(t_1,t_2)) \\
&=& \sum_{v\in V(\Gamma)} e^{w_{D_\beta}^v  - w_{D_\alpha}^v}
N_{Y_\alpha^v, Y_\beta^v}(w_1^v, w_2^v)
\end{eqnarray*}
where $w_1^v=w_1^v(\ep_1,\ep_2)$, $w_2^v=w_2^v(\ep_1,\ep_2)$,
  $t_1=e^{\ep_1}, t_2=e^{\ep_2}$.
\end{pro}
\begin{proof} 
\begin{equation}\label{eqn:M}
\begin{aligned}
M_{\alpha,\beta}(t_1,t_2)=& \ch_{T_t}\Ext^*(\cO_X(D_\alpha),\cO_X(D_\beta-\linf))\\
& -\ch_{T_t}\Ext^*(I_\alpha(D_\alpha),I_\beta(D_\beta-\linf)) .
\end{aligned}
\end{equation}
We will compute the two terms on the right hand side of \eqref{eqn:M}
using the method in \cite[Section 4]{MNOP1}.
\begin{eqnarray*}
&& \Ext^*(I_\alpha(D_\alpha), I_\beta(D_\beta  - \linf))\\\\
&=& \sum_{i,j=0}^2 (-1)^{i+j} H^i(X,\ext^j(I_\alpha(D_\alpha), I_\beta(D_\beta-\linf))\\
&=& \sum_{i,j=0}^2 (-1)^{i+j} \fC^i(X,\ext^j(I_\alpha(D_\alpha), I_\beta(D_\beta-\linf))
\end{eqnarray*}
where $\fC^i$ denote the \v{C}ech cochain groups. More explicitly, let
$$
\{ p_a\mid a=1,\ldots,\chi(X)\}
$$
be the $T_t$-fixed points in $X$, where $\chi(X)$ is the Euler characteristic
of $X$.
Let $U_a$ be the $\bC^2$ coordinate chart with origin at $p_a$, and let
$U_{ab }=U_a\cap U_b$, etc.

\begin{eqnarray*}
&& \sum_{i,j=0}^2 (-1)^{i+j} 
\fC^i\left(X,\ext^j(I_\alpha(D_\alpha), I_\beta(D_\beta-\linf))\right)\\
&=& \bigoplus_{a}\sum_{j=0}^2 (-1)^j 
\Gamma\left(U_a,\ext^j(I_\alpha(D_\alpha), I_\beta(D_\beta-\linf))\right)\\
&&  -\bigoplus_{a,b}\sum_{j=0}^2 (-1)^j 
\Gamma\left(U_{ab},\ext^j(I_\alpha(D_\alpha), I_\beta(D_\beta-\linf))\right)\\
&&  +\bigoplus_{a,b,c}\sum_{j=0}^2 (-1)^j 
\Gamma\left(U_{abc},\ext^j(I_\alpha(D_\alpha), I_\beta(D_\beta-\linf))\right)\\
&&\ldots 
\end{eqnarray*}
Note that $I_\alpha|_{U_{a_1 ... a_i}} =\cO_X|_{U_{a_1 ... a_i}}$ unless
$i=1$ and $p_{a_1}\in X_0$, so 
\begin{eqnarray*}
&& \Ext^*(\cO_X(D_\alpha), \cO_X(D_\beta  - \linf))
-\Ext^*(I_\alpha(D_\alpha), I_\alpha(D_\beta  - \linf))\\
&=&  \bigoplus_{v\in V(\Gamma)}\sum_{j=0}^2 (-1)^j 
\Gamma\left(U_v,\ext^j(\cO_X(D_\alpha), \cO_X(D_\beta))\right)\\
&& - \bigoplus_{v\in V(\Gamma)}\sum_{j=0}^2 (-1)^j 
\Gamma\left(U_v,\ext^j(I_\alpha(D_\alpha), I_\beta(D_\beta))\right)
\end{eqnarray*}
where $U_v$ is the $\bC^2$ chart centered at $p_v$.

Given a partition (Young diagram) $Y$ and a box $x\in Y$,
let $a'(x)$ and $l'(x)$ be the arm-colength 
and leg-colength of $x$, respectively (see e.g. \cite[Section 3.1]{NY}).
Given a partition $Y$,  we define
$$
Q_Y(s_1,s_2)=\sum_{x\in Y} s_1^{l'(x)} s_2^{a'(x)}.
$$
We have
\begin{eqnarray*}
&& \ch_{\tT}\sum_{j=0}^2 (-1)^j 
\Gamma\left(U_v,\ext^j(\cO_X(D_\alpha), \cO_X(D_\beta))\right)\\
&& -\ch_{\tT}\sum_{j=0}^2 (-1)^j 
\Gamma\left(U_v,\ext^j(I_\alpha(D_\alpha), I_\beta(D_\beta))\right)\\
&=& \chi^v_{D_\beta} \bigl(\chi^v_{D_\alpha}\bigr)^{-1}
\Bigl( Q_{Y_\alpha^v}(\chi_1^v,\chi_2^v)\chi_1^v  \chi_2^v
 +Q_{Y^v_\beta}( (\chi_1^v)^{-1}, (\chi_2^v)^{-1}) \\
&& -Q_{Y^v_\alpha}(\chi_1^v,\chi_2^v) 
Q_{Y^v_\beta}((\chi_1^v)^{-1},(\chi_2^v)^{-1})(1-\chi_1^v)(1-\chi_2^v)\Bigr)\\
&=& \chi^v_{D_\alpha}(\chi^v_{D_\alpha})^{-1}
M_{Y^v_\alpha,Y^v_\beta}(\chi_1^v,\chi_2^v)
\end{eqnarray*}
where
\begin{eqnarray*}
M_{S,T}(t_1,t_2)&=&
 Q_{S}(t_1,t_2)t_1 t_2 
+Q_{T}(t_1^{-1},t_2^{-1})\\
&&  -Q_{S}(t_1,t_2) Q_{T}(t_1^{-1},t_2^{-1})(1-t_1)(1-t_2).
\end{eqnarray*}
We now compare our expression of $M_{Y^v_\alpha,Y^v_\beta}(t_1,t_2)$ with
the notation in the proof of \cite[Theorem 2.11]{NY1}. The correspondence is
\begin{eqnarray*}
t_1 t_2 \Hom(V_\alpha,W_\beta)&=& Q_{Y_\alpha^v}(t_1,t_2)t_1 t_2 \\
\Hom(W_\alpha,V_\beta)&=& Q_{Y^v_\beta}(t_1^{-1},t_2^{-1})\\
(t_1+t_2-1-t_1t_2)\Hom(V_\alpha,V_\beta)&=&
-Q_{Y^v_\alpha}(t_1,t_2) Q_{Y^v_\beta}(t_1^{-1},t_2^{-1})(1-t_1)(1-t_2)
.\end{eqnarray*}
So $M_{S,T}(t_1,t_2)$ can be rewritten as \eqref{eqn:MST}.
\end{proof}

\subsubsection{$L_{\alpha,\beta}(t_1,t_2)$}
\begin{lm} \label{Lvanish}
If $D_\alpha=D_\beta$ then
$L_{\alpha,\beta}(t_1,t_2)=0$. In particular,
$L_{\alpha,\alpha}(t_1,t_2)=0$.
\end{lm}
\begin{proof}
$$
L_{\alpha,\beta}(t_1,t_2)=-\chi_{T_t}(X,\cO_X(-\linf))
$$
which can be identified with the tangent space of
$\fM_{1,0,0}(X,\linf)$ at the trivial line bundle
$\cO_X$. By Proposition \ref{thm:vanish},
$$
H^0(X,\cO_X(-\linf))=H^2(X,\cO_X(-\linf))=0.
$$
By Corollay \ref{thm:dim} (here $r=1,d=0, n=0$), $H^1(X,\cO_X(-\linf))=0$.
\end{proof}
By Proposition \ref{thm:vanishV} and Corollary \ref{thm:rank},
we have 
\begin{lm}\label{Dvanish} Suppose that $D\cdot \linf =0$. Then
$$
H^0(X,\cO_X(D-\linf))=  H^2(X,\cO_X(D -\linf))=0,
$$
and
$$
\dim_{\bC}H^1(X,\cO_X(D-\linf)) =
-\frac{1}{2}\left(D^2 + c_1(X)\cdot D\right).
$$
\end{lm}
In particular, for any $D$ such that $D\cdot\linf=0$ we have
$$
D^2 = - \dim_\bC H^1 (X,\cO_X(D-\linf)) - \dim_\bC H^1(X,\cO_X(-D-\linf)) \leq 0.
$$

\begin{nt}\label{nt:weights}
Let $q_0,$$q_1$ be  the two $T_t$ fixed points on $\linf.$
Let $w$ (resp. $u$) $\in \mathrm{Hom}(T,\bC^*)$ be the {\em tangent weight} (resp. {\em normal weight}) at $q_0$, i.e., the weight of the $T_t$-action on
$T_{q_0} \linf$ (resp. $(N_{\linf/X})_{q_0}$). Then the 
tangent weight (resp. normal weight) at $q_1$, i.e.,
the weight of the $T_t$-action on $T_{q_1}\linf$
(resp. $(N_{\linf/X})_{q_1}$), must be given by
$-w$ (resp. $u-kw$), where $k=\linf\cdot\linf >0$.
\end{nt}

Note that the normal weights at $q_0$ and $q_1$
are the restrictions of the equivariant 
first Chern class $(c_1)_{T_t}(\cO_X(\linf))$
to the $T_t$ fixed points $q_0$ and $q_1$, respectively:
$$
(c_1)_{T_t}(\cO_X(\linf))\Bigr|_{q_0} = u,\quad
(c_1)_{T_t}(\cO_X(\linf))\Bigr|_{q_1} = u-kw.
$$

\begin{pro}[edge contribution to the tangent bundle] \label{thm:T-L}
$$
L_{\alpha,\beta}(t_1,t_2)
= \left(\sum_{v\in V(\Gamma)} 
\frac{ -e^{w_{D_\beta}^v-w_{D_\alpha}^v }  }{(1-e^{-w_1^v})(1-e^{-w_2^v}) } \right)
+\frac{1}{ (1-e^{-w})(1-e^u)}+\frac{1}{(1-e^w)(1-e^{u-kw})}. 
$$
\end{pro}
\begin{proof} Recall that $L_{\alpha,\beta}(t_1,t_2)= 
-\chi_{T_t}(X,\cO_X(D_\beta-D_\alpha-\linf) )$. By 
Grothendieck-Riemann-Roch,
\begin{eqnarray*}
&& \chi_{T_t}(X,\cO_X (D_\beta-D_\alpha-\linf))\\
&=& \int_X \td_{T_t}(T_X)\ch_{T_t}(\cO_X(D_\beta-D_\alpha-\linf)) \\
&=&  \left(\sum_{v\in V(\Gamma)} 
\frac{ e^{w_{D_\beta}^v-w_{D_\alpha}^v }  }{(1-e^{-w_1^v})(1-e^{-w_2^v}) } \right)
+\frac{e^{-u}}{ (1-e^{-w})(1-e^{-u})}+\frac{e^{-u+kw}}{(1-e^w)(1-e^{-u+kw})}. 
\end{eqnarray*}
\end{proof}

\begin{ex}\label{FkL}
Let $X=\bF_k$, $\lz$, $\linf$ be as in Example \ref{Fk}, with
the following $T_t$-action:

\bigskip

\begin{tabular}{c|c|c|c|c|c|c|c}
$T_{p_1}\lz$ & $(N_{\lz/X})_{p_1}$ & $T_{p_2}\lz$ & $(N_{\lz/X})_{p_2}$ &
$T_{p_3}\linf$ & $(N_{\linf/X})_{p_3}$ & $T_{p_4}\linf$ & $(N_{\linf/X})_{p_4}$ \\ 
\hline
$\ep_1$ & $\ep_2$ &  $-\ep_1$ & $\ep_2+ k\ep_1$ & $-\ep_1$ & $-\ep_2-k\ep_1$ & $\ep_1$ & $-\ep_2$
\end{tabular}

\bigskip

Hence, here $w=\ep_1$ and $u=-\ep_2$, and we have $D_\alpha= d_\alpha \ell_0$ for some $d_\alpha\in\bZ$. Then
\begin{eqnarray*}
L_{\alpha,\beta}(t_1,t_2)
&=& \frac{ -e^{ (d_\beta-d_\alpha)\ep_2}  }{(1-e^{-\ep_1})(1-e^{-\ep_2}) } 
+\frac{-e^{ (d_\beta-d_\alpha)(\ep_2+k\ep_1) }}{(1-e^{\ep_1})(1-e^{-\ep_2-k\ep_1}) }\\
&& +\frac{1}{ (1-e^{-\ep_1})(1-e^{-\ep_2})}+\frac{1}{(1-e^{\ep_1})(1-e^{-\ep_2-k\ep_1})}\\
&=& \frac{ 1-t_2^{ d_\beta-d_\alpha}  }{(1-t_1^{-1})(1-t_2^{-1}) } 
+\frac{1-(t_1^k t_2)^{ d_\beta-d_\alpha}}{(1-t_1)(1-t_1^{-k}t_2^{-1}) }
,\end{eqnarray*} and 
we have
$$
L_{\alpha,\beta}(t_1,t_2)=\begin{cases}
\displaystyle{\sum_{j=0}^{d_\alpha-d_\beta-1} \sum_{i=0}^{kj} t_1^{-i} t_2^{-j}  }& \textup{if }d_\alpha>d_\beta,\\
\displaystyle{\sum_{j=1}^{d_\beta-d_\alpha}\sum_{i=1}^{kj-1} t_1^i t_2^j }& \textup{if }d_\alpha<d_\beta,\\
0 & \textup{if }d_\alpha=d_\beta.
\end{cases} 
$$
\end{ex}

\subsection{The natural bundle: fundamental representation}\label{sec:Vinstanton}

Let $(E,\Phi)\in \MX$ be a fixed point of the  $\tT$-action corresponding to
$\DY=(D_1,\vY_1,\ldots,D_r, \vY_r)$. 
We want to compute
$$
\ch_{\tT}V_{(E,\Phi)}  = \ch_{\tT} H^1(X,E(-\linf))
=-\chi_{\tT} (X,E(-\linf)).
$$ 
We have
$$
E= I_1(D_1)\oplus \cdots \oplus I_r(D_r),
$$
so 
$$
-\chi_{\tT}(X,E(-\linf)) 
= -\sum_\beta\chi_{\tT}(X,I_\beta(D_\beta-\linf))
=-\sum_\beta e^{a_\beta} \chi_{T_t}(X,I_\beta(D_\beta-\linf)).
$$
Let
\begin{eqnarray*}
L_\beta(t_1,t_2)&=& -\chi_{T_t}(X,\cO_X(D_\beta-\linf))\\
M_\beta(t_1,t_2)&=& \chi_{T_t}(X,\cO_X(D_\beta-\linf))
-\chi_{T_t}(X,I_\beta(D_\beta-\linf)). 
\end{eqnarray*}
Then
\begin{equation}\label{eqn:V-ML}
\ch_{\tT} V_{(E,\Phi)} =\sum_{\beta=1}^r e^{a_\beta}
\left( M_\beta(t_1,t_2) + L_\beta(t_1,t_2)\right).
\end{equation}
So it remains to compute $M_\beta(t_1,t_2)$ and $L_\beta(t_1,t_2)$.

Let $w_{D_\alpha}^v, w_1^v, w_2^v$ be defined as in
Section \ref{sec:M-T}.
Given a partition $S$, let
\begin{equation}
M_S(t_1,t_2)=\sum_{s\in S} t_1^{-l'(s)} t_2^{-a'(s)}
\end{equation}
\begin{equation}
N_S(\ep_1,\ep_2) \define M_S(e^{\ep_1}, e^{\ep_2})
 = \sum_{s\in S} e^{ -l'(s)\ep_1 - a'(s)\ep_2}.
\end{equation}

\begin{pro}[vertex contribution to the natural bundle]\label{thm:V-M}
\begin{eqnarray*}
M_\beta(t_1,t_2)&=&\sum_{v\in V(\Gamma)}
\chi_{D_\beta}^v(t_1,t_2) M_{Y_\beta^v}(\chi_1^v(t_1,t_2), \chi_2^v(t_1,t_2)) \\
&=& \sum_{v\in V(\Gamma)} e^{w_{D_\beta}^v}
N_{Y_\beta^v}(w_1^v, w_2^v).
\end{eqnarray*}
\end{pro}
\begin{proof} 
Let $D_\alpha=0$ in Proposition \ref{thm:T-M}.
\end{proof}

\begin{pro}[edge contribution to the natural bundle]\label{thm:V-L}
$$
L_{\beta}(t_1,t_2)
= \left(\sum_{v\in V(\Gamma)} 
\frac{ -e^{w_{D_\beta}^v}  }{(1-e^{-w_1^v})(1-e^{-w_2^v}) } \right)
+\frac{1}{ (1-e^{-w})(1-e^u)}+\frac{1}{(1-e^w)(1-e^{u-kw})}. 
$$
\end{pro}
\begin{proof} 
Let $D_\alpha=0$ in Proposition \ref{thm:T-L}.
\end{proof}

\begin{ex}
Let $X=\bF_k$, $\lz$, $\linf$ be as in Example \ref{Fk}, with
the $T_t$-action as in Example \ref{FkL}.
Then
$$
L_{\beta}(t_1,t_2)=\begin{cases}
\displaystyle{\sum_{j=0}^{-d_\beta-1} \sum_{i=0}^{kj} t_1^{-i} t_2^{-j}  }& \textup{if }d_\beta<0,\\
\displaystyle{\sum_{j=1}^{d_\beta}\sum_{i=1}^{kj-1} t_1^i t_2^j }& \textup{if }d_\beta>0,\\
0 & \textup{if }d_\beta=0.
\end{cases} 
$$
\end{ex}

\subsection{Formula for instanton partition functions}

Given $\vY=(Y_1,\ldots,Y_r)$, where each $Y_\alpha$ is a Young diagram,
and a multiplicative class $A$ associated to $f(x)$, define
\begin{equation}\label{eqn:m}
\begin{aligned}
m_{A,\alpha,\beta}^\vY(\ep_1,\ep_2,\va)
\define &\prod_{s\in Y_\alpha}f(a_\beta-a_\alpha -l_{Y_\beta}(s) \ep_1 +(a_{Y_\alpha}(s)+1)\ep_2)\\
&\cdot \prod_{t\in Y_\beta}f(a_\beta-a_\alpha+(l_{Y_\alpha}(t)+1)\ep_1 -a_{Y_\beta}(t)\ep_2),
\end{aligned}
\end{equation}
\begin{equation} \label{eqn:m-beta}
m_{A,\beta}^\vY(\ep_1,\ep_2,\va)\define
\prod_{t\in Y^\beta} f(a_\beta -l'_{Y_\beta}(t)\ep_1 -a'_{Y_\beta}(t)\ep_2).
\end{equation}
In particular,
\begin{equation}\label{eqn:m-ctop}
\begin{aligned}
m_{c_{top},\alpha,\beta}^\vY(\ep_1,\ep_2,\va)
= &\prod_{s\in Y_\alpha}(a_\beta-a_\alpha -l_{Y_\beta}(s) \ep_1 +(a_{Y_\alpha}(s)+1)\ep_2)\\
&\cdot \prod_{t\in Y_\beta}(a_\beta-a_\alpha+(l_{Y_\alpha}(t)+1)\ep_1 -a_{Y_\beta}(t)\ep_2).
\end{aligned}
\end{equation}

Let $Z_{\bC^2,A,B}^\inst = Z^\inst_{\bC^2,A,B,0}$, and let
$|\vY|=\sum_{\alpha=1}^r |Y_\alpha|$. In this case, 
all $D_\beta=0$, so the leg contribution is zero (see Lemma \ref{Lvanish}, 
Lemma \ref{Dvanish}):
$$
L_{\alpha,\beta}=0, \quad  L_\beta=0.
$$
By \eqref{eqn:T-ML}, Proposition \ref{thm:T-M}, \eqref{eqn:V-ML}, 
Proposition \ref{thm:V-M}, and 
above definitions \eqref{eqn:m}, \eqref{eqn:m-beta}, \eqref{eqn:m-ctop}, we have:
\begin{pro}[instanton partition functions for $\bC^2$]\label{thm:CtwoAB}
$$
Z^\inst_{\bC^2,A,B}(\ep_1,\ep_2,\va;\Lambda)=\sum_{\vY}
\Lambda^{2r|\vY|}\prod_{\alpha,\beta}
\frac{m_{A,\alpha,\beta}^\vY(\ep_1,\ep_2,\va) }
{m_{c_{\mathrm{top}},\alpha,\beta}^\vY(\ep_1,\ep_2,\va)} 
\prod_{\beta=1}^r m_{B,\beta}^\vY(\ep_1,\ep_2,\va) 
$$
\end{pro}

Given $\vD=(D_1,\ldots,D_r)$, where each $D_\alpha \in \oplus_{e\in E(\Gamma)} \bZ \ell_e \cong H_2(X_0;\bZ)$,
and a multiplicative class $A$, define
\begin{equation}\label{eqn:l}
l_{A,\alpha,\beta}^{\vD}(\ep_1,\ep_2,\va)= A_\tT H^1(X,\cO_X(D_\beta-D_\alpha-\linf)). 
\end{equation}
Then $l_{A,\alpha,\beta}^\vD(\ep_1,\ep_2;\va)=1$ if $D_\alpha=D_\beta$.
In particular, $l_{A,\alpha,\alpha}^\vD(\ep_1,\ep_2;\va)=1$.
Let
\begin{equation}
l_{A,\beta}^{\vD}(\ep_1,\ep_2,\va)= A_\tT H^1(X,\cO_X(D_\beta-\linf)).
\end{equation}
Let
$$
|\vD|^2 =-\frac{1}{2}\sum_{\alpha\neq \beta}(D_\alpha-D_\beta)^2 \geq 0.
$$

By Equations \eqref{eqn:T-ML}, \eqref{eqn:V-ML} and
Propositions \ref{thm:T-M}, \ref{thm:T-L}, \ref{thm:V-M}, \ref{thm:V-L},
\ref{thm:CtwoAB}, we have the following analogue
of the ``master formula'' in \cite[Section 6]{Ne3}.
 
\begin{pro}[master formula for instanton partition functions]\label{master}
\begin{eqnarray*}
Z^\inst_{X_0,A,B,d}(\ep_1,\ep_2,\va;\Lambda) 
&=&\sum_{\sum D_\alpha =d } \Lambda^{|\vD|^2}\prod_{\alpha\neq \beta}
\frac{l_{A,\alpha,\beta}^\vD(\ep_1,\ep_2,\va)}
{l_{c_{top},\alpha,\beta}^\vD(\ep_1,\ep_2,\va)}
 \prod_{\beta=1}^r l_{B,\beta}^\vD(\ep_1,\ep_2,\va)\\
&& \quad\cdot \prod_{v\in V(\Gamma)} Z_{\bC^2,A,B}^\inst(w_1^v, w_2^v, \va+ \vD^v;\Lambda)
\end{eqnarray*}
where $\vD^v=(w_{D_1}^v,\ldots, w_{D_r}^v)$.
\begin{eqnarray*}
Z^\inst_{X_0,A,B}(\ep_1,\ep_2,\va;\Lambda,Q) 
&=&\sum_{D_\alpha\in H^2_c(X;\bZ)} Q^{\sum_\alpha D_\alpha} \Lambda^{|\vD|^2}\prod_{\alpha\neq \beta}
\frac{l_{A,\alpha,\beta}^\vD(\ep_1,\ep_2,\va)}
{l_{c_{top},\alpha,\beta}^\vD(\ep_1,\ep_2,\va)}\\
&&\cdot \prod_{\beta=1}^r l_{B,\beta}^\vD(\ep_1,\ep_2,\va)
\cdot \prod_{v\in V(\Gamma)} Z_{\bC^2,A,B}^\inst(w_1^v, w_2^v, \va+ \vD^v;\Lambda)
\end{eqnarray*}
\end{pro}
In the rank 1 case, $Z^\inst_{X_0,A,B}$ does not depend on $\va$.
\begin{cor}[rank 1, $B=1$ case]\label{rank-one}
\begin{eqnarray*}
Z^\inst_{X_0,A,B=1,d}(\ep_1,\ep_2;\Lambda) 
&=&\prod_{v\in V(\Gamma)} Z_{\bC^2,A,B=1}^\inst(w_1^v, w_2^v;\Lambda)\\
Z^\inst_{X_0,A,B=1}(\ep_1,\ep_2;\Lambda,Q) 
&=&\sum_{d \in H^2_c(X;\bZ)} Q^{d} 
\prod_{v\in V(\Gamma)} Z_{\bC^2,A,B=1}^\inst(w_1^v, w_2^v;\Lambda)
\end{eqnarray*}
\end{cor}

Note that Corollary \ref{rank-one} is applicable to the
following cases: 4d pure gauge theory (Section \ref{sec:four-pure}), 
4d gauge theory with one adjoint matter hypermultiplet (Section
\ref{sec:four-adjoint}), 5d gauge theory compactified 
on a circle (Section \ref{sec:five}), 
Hirzebruch genus (Section \ref{sec:Hirzebruch}), 
elliptic genus (Section \ref{sec:elliptic}).

\subsection{Nekrasov  conjecture for $\bC^2$: instanton part}

\begin{df}[instanton prepotential for $\bC^2$]
Define 
$$
\cF^\inst_{\bC^2,A,B}(\ep_1,\ep_2,\va;\Lambda)\define
-\ep_1\ep_2 \log Z^\inst_{\bC^2,A,B}(\ep_1,\ep_2,\va;\Lambda).
$$
\end{df}

There are several versions of Nekrasov conjecture
which correspond to the following special cases:
\begin{enumerate}
\item 4d pure gauge theory 
(see Section \ref{sec:four-pure}):
$$
\cF^\inst_{\bC^2}(\ep_1,\ep_2,\va;\Lambda)=
\cF^\inst_{\bC^2, A=1, B=1}(\ep_1,\ep_2,\va;\Lambda).
$$
\item 4d gauge theory with $N_f$ fundamental matter hypermultiplets
(see Section \ref{sec:four-fundamental}):
$$
\cF^\inst_{\bC^2}(\ep_1,\ep_2,\va,\vm;\Lambda)=
\cF^\inst_{\bC^2, A=1, B=E_{\vm} }(\ep_1,\ep_2,\va;\Lambda).
$$

\item 4d gauge theory with one adjoint matter hypermultiplet
(see Section \ref{sec:four-adjoint}):
$$
\cF^\inst_{\bC^2}(\ep_1,\ep_2,\va,m;\Lambda)=
\cF^\inst_{\bC^2, A=E_m, B=1}(\ep_1,\ep_2,\va,m;\Lambda).
$$

\item 5d gauge theory compactified on a circle of circumference 
$\beta$ (see Section \ref{sec:five}):
$$
\cF^\inst_{\bC^2}(\ep_1,\ep_2,\va;\Lambda,\beta)=
\cF^\inst_{\bC^2, A=\widehat{A}_\beta, B=1}(\ep_1,\ep_2,\va,m;\Lambda).
$$
\end{enumerate}

The above definitions of $\cF^\inst_{\bC^2}$ are the same as those
in \cite{NO}; the definition in case (1) above is the negative
of the definition in \cite{NY, NY1}.

In Theorem \ref{thm:Nekrasov-Ctwo} below, 
we summarize the various versions of the Nekrasov conjecture
proved by Nakajima-Yoshioka \cite{NY1, NY2}, 
Nekrasov-Okounkov \cite{NO}, Braverman-Etingof \cite{Br, BE},
G\"{o}ttsche-Nakajima-Yoshioka \cite{GNY2}. We refer
to Appendix \ref{sec:SW} for the definitions
of the corresponding versions of the Seiberg-Witten prepotential
in Theorem \ref{thm:Nekrasov-Ctwo}.

\begin{thm}[Nekrasov conjecture for $\bC^2$: instanton part]
\label{thm:Nekrasov-Ctwo}
$ $\\
\begin{enumerate}
\item {\em 4d pure gauge theory \cite{NO, NY1, BE}:}
\begin{enumerate}
\item $\cF^\inst_{\bC^2}(\ep_1,\ep_2,\va;\Lambda)$ is analytic
in $\ep_1,\ep_2$ near $\ep_1=\ep_2=0$.
\item $\displaystyle{\lim_{\ep_1,\ep_2\to 0} \cF^\inst_{\bC^2}(\ep_1,\ep_2,\va;\Lambda)}
= \cF_0^\inst(\va,\Lambda)$, where
$\cF_0^\inst(\va,\Lambda)$ is the instanton part of the Seiberg-Witten prepotential
of 4d pure gauge theory.
\end{enumerate}

\item {\em 4d gauge theory with $N_f$ fundamental matter 
hypermultiplets \cite{NO}:}
\begin{enumerate}
\item $\cF^\inst_{\bC^2}(\ep_1,\ep_2,\va,\vm;\Lambda)$ is analytic
in $\ep_1,\ep_2$ near $\ep_1=\ep_2=0$.
\item $\displaystyle{\lim_{\ep_1,\ep_2\to 0} \cF^\inst_{\bC^2}(\ep_1,\ep_2,\va,\vm;\Lambda)}
= \cF_0^\inst(\va,\vm,\Lambda)$, where
$\cF_0^\inst(\va,\vm,\Lambda)$ is the instanton part of the Seiberg-Witten prepotential
of 4d gauge theory with $N_f$ fundamental matter hypermultiplets.
\end{enumerate}

\item {\em 4d gauge theory with one adjoint matter hypermultiplet \cite{NO}:}
\begin{enumerate}
\item $\cF^\inst_{\bC^2}(\ep_1,\ep_2,\va,m;\Lambda)$ is analytic
in $\ep_1,\ep_2$ near $\ep_1=\ep_2=0$.
\item $\displaystyle{\lim_{\ep_1,\ep_2\to 0} \cF^\inst_{\bC^2}(\ep_1,\ep_2,\va,m;\Lambda)}
= \cF_0^\inst(\va,m,\Lambda)$, where
$\cF_0^\inst(\va,m,\Lambda)$ is the instanton part of the Seiberg-Witten prepotential
of 4d gauge theory with one adjoint matter hypermultiplet.
\end{enumerate}

\item {\em 5d gauge theory compactified on a circle of circumference 
$\beta$ \cite{NO, NY2, GNY2}:}
\begin{enumerate}
\item $\cF^\inst_{\bC^2}(\ep_1,\ep_2,\va;\Lambda,\beta)$ is analytic
in $\ep_1,\ep_2$ near $\ep_1=\ep_2=0$.
\item $\displaystyle{\lim_{\ep_1,\ep_2\to 0} \cF^\inst_{\bC^2}(\ep_1,\ep_2,\va;\Lambda,\beta)}
= \cF_0^\inst(\va,\Lambda,\beta)$, where
$\cF_0^\inst(\va,\Lambda,\beta)$ is the instanton part of the Seiberg-Witten prepotential
of 5d gauge theory compactified on a circle of circumference $\beta$.
\end{enumerate}
\end{enumerate}
\end{thm}

\subsection{Nekrasov conjecture for toric surfaces: instanton part}
The expression of the master formula (Proposition  \ref{master}) contains  two parts. 
\begin{itemize}
\item Leg contribution:
$$
\prod_{\alpha\neq\beta}
\frac{l^\vD_{A,\alpha,\beta}(\ep_1,\ep_2,\va)}
{l^{\vD}_{c_{top},\alpha,\beta}(\ep_1,\ep_2,\va) }
\prod_{\beta=1}^r l_\beta^\vD(\ep_1,\ep_2,\va)
$$
is analytic in $\ep_1, \ep_2$ near $\ep_1,\ep_2=0$, and
\begin{eqnarray*}
&& \lim_{\ep_1,\ep_2\to 0}
\prod_{\alpha\neq\beta}
\frac{l^\vD_{A,\alpha,\beta}(\ep_1,\ep_2,\va)}
{l^{\vD}_{c_{top},\alpha,\beta}(\ep_1,\ep_2,\va) }
\prod_{\beta=1}^r l_\beta^\vD(\ep_1,\ep_2,\va)\\
&& =\prod_{\alpha\neq \beta} 
\left(\frac{f(a_\beta-a_\alpha)}{a_\beta-a_\alpha}\right)
^{-\frac{1}{2}((D_\beta-D_\alpha)^2+c_1(X)(D_\beta-D_\alpha) )}
\prod_{\beta=1}^r g(a_\beta)^{-\frac{1}{2}
( D_\beta^2 + c_1(X)\cdot D_\beta)}.
\end{eqnarray*}
\item Vertex contribution:
$$
\prod_{v\in V(\Gamma)} Z^\inst_{\bC^2,A,B}(w_1^v, w_2^v, \va+\vD^v;\Lambda)
=\exp\left( - \sum_{v\in V(\Gamma) }
\frac{\cF^\inst_{\bC^2,A,B}(w_1^v,w_2^v,\va+\vD^v;\Lambda)}{w_1^v w_2^v}  \right).
$$
\end{itemize}

\begin{df}
Given $\vD=(D_1,\ldots,D_r)$, where
each  $D_\alpha \in \displaystyle{\bigoplus_{e\in E(\Gamma)} \bZ \ell_e} = H_2(X_0;\bZ)$, define
\begin{eqnarray*}
 \cF^\inst_{X_0,A,B,\vD}(\ep_1,\ep_2,\va;\Lambda)
&=&\sum_{v\in V(\Gamma)}\frac{\cF^\inst_{\bC^2,A,B}(w_1^v,w_2^v,\va+\vD^v;\Lambda)}{w_1^v w_2^v}\\
&& +\frac{\cF^\inst_{\bC^2,A,B}(w,u,\va;\Lambda) }{wu}
+\frac{\cF^\inst_{\bC^2,A,B}(-w,u-kw,\va;\Lambda)}{-w(u-kw)}.
\end{eqnarray*}
\end{df}

\begin{lm}\label{thm:ABD}
Assume that $\cF^\inst_{\bC^2, A,B}(\ep_1,\ep_2,\va;\Lambda)$
is analytic in $\ep_1,\ep_2$ near $\ep_1=\ep_2=0$. Then
$\cF^\inst_{X_0,A,B,\vD}(\ep_1,\ep_2,\va;\Lambda)$ is
analytic in $\ep_1,\ep_2$ near $\ep_1=\ep_2=0$ for
all $\vD$.
\end{lm}
\begin{proof}
$\cF_{\bC^2,A,B}^\inst(\ep_1,\ep_2,\va;\Lambda)$ is symmetric in $\ep_1, \ep_2$,
so it is a function of $s_1=\ep_1+ \ep_2$, $s_2= \ep_1\ep_2$, $\va$, and
$\Lambda$. For fixed $\va,\Lambda$, let
$$
g_{A,B}(s_1, s_2, d_1,\ldots, d_r,\va;\Lambda)=
\cF_{\bC^2,A,B}^\inst(\ep_1,\ep_2,a_1+d_1,\ldots, a_r + d_r;\Lambda).
$$
Then $g_{A,B}(s_1,s_2,d_1,\ldots,d_r,\va;\Lambda)$ is analytic in $s_1, s_2, d_1,\ldots, d_r$ near
$s_1=s_2=d_1=\cdots= d_r=0$, so it has a power series expansion.
Therefore
\begin{equation}\label{eqn:cTt}
\begin{aligned}
& I_{A,B,\vD}(\ep_1,\ep_2;\va,\Lambda)\\
\define & \int_X g_{A,B}\Bigl((c_1)_{T_t}(T_X), (c_2)_{T_t}(T_x), 
(c_1)_{T_t}(\cO_X(D_1)),\ldots, (c_1)_{T_t}(\cO_X(D_r)),\va;\Lambda\Bigr)
\end{aligned}
\end{equation}
is analytic in $\ep_1, \ep_2$ near $\ep_1 =\ep_2 =0$, and 
\begin{equation}\label{eqn:cTt-limit}
\lim_{\ep_1,\ep_2\to 0} I_{A,B,\vD}(\ep_1,\ep_2;\va,\Lambda)
=\int_X g_{A,B}\Bigl(c_1(T_X), c_2(T_x), c_1(\cO_X(D_1)),\ldots, c_1(\cO_X(D_r)),\va;\Lambda\Bigr).
\end{equation}
The integral $I_{A,B,\vD}(\ep_1,\ep_2,\va;\Lambda)$ is computed by the localization formula as follows: 
\begin{eqnarray*}
I_{A,B,\vD}(\ep_1,\ep_2,\va;\Lambda)&=&
\sum_{v\in V(\Gamma)}\frac{\cF^\inst_{\bC^2,A,B}(w_1^v,w_2^v,\va+\vD^v;\Lambda)}{w_1^v w_2^v}\\
&& +\frac{\cF^\inst_{\bC^2,A,B}(w,u,\va;\Lambda) }{wu}
+\frac{\cF^\inst_{\bC^2,A,B}(-w,u-kw,\va;\Lambda)}{-w(u-kw)}
\end{eqnarray*}
\end{proof}

\begin{df}
Assume that $\cF^\inst_{\bC^2, A,B}(\ep_1,\ep_2,\va;\Lambda)$
is analytic in $\ep_1,\ep_2$ near $\ep_1=\ep_2=0$. Define
$$
F_{X_0,A,B,\vD}(\va;\Lambda)\define
\lim_{\ep_1,\ep_2\to 0} \cF^\inst_{X_0,A,B,\vD}(\ep_1,\ep_2,\va;\Lambda).
$$
\end{df}

\begin{lm} \label{thm:logloglog}
If $\cF^\inst_{\bC^2,A,B}(\ep_1,\ep_2,\va;\Lambda)$ is analytic
in $\ep_1, \ep_2$ near $\ep_1=\ep_2=0$, then 
$$ 
\log\Bigl(Z^\inst_{X_0,A,B,d}(\ep_1,\ep_2;\va;\Lambda)
Z^\inst_{\bC^2,A,B}(w,u,\va;\Lambda)
Z^\inst_{\bC^2,A,B}(-w,u-kw,\va;\Lambda)\Bigr)
$$
is analytic in $\ep_1,\ep_2$ near
$\ep_1=\ep_2=0$.
\end{lm}
\begin{proof} We have
\begin{eqnarray*}
&& Z^\inst_{X_0,A,B,d}(\ep_1,\ep_2;\va;\Lambda)
Z^\inst_{\bC^2,A,B}(w,u,\va;\Lambda)
Z^\inst_{\bC^2,A,B}(-w,u-kw,\va;\Lambda)\\
&=& \sum_{\sum D_\alpha =d}\Lambda^{|\vD|^2} 
h_\vD(\ep_1,\ep_2,\va;\Lambda)
\end{eqnarray*}
where
\begin{eqnarray*}
h_\vD(\ep_1,\ep_2,\va;\Lambda)= \prod_{\alpha\neq\beta}
\frac{l^\vD_{A,\alpha,\beta}(\ep_1,\ep_2,\va)}
{l^{\vD}_{c_{top},\alpha,\beta}(\ep_1,\ep_2,\va) }
\prod_{\beta=1}^r l_\beta^\vD(\ep_1,\ep_2,\va)\exp\left(
-\cF^\inst_{X,A,B,\vD}(\ep_1,\ep_2,\va;\Lambda)\right).
\end{eqnarray*}
$h_\vD(\ep_1,\ep_2,\va;\Lambda)$ is analytic in
$\ep_1, \ep_2$ near $\ep_1=\ep_2=0$, and
\begin{eqnarray*}
 \lim_{\ep_1,\ep_2\to 0} h_\vD(\ep_1,\ep_2,\va;\Lambda)&=&
\prod_{\alpha\neq \beta} 
\left(\frac{f(a_\beta-a_\alpha)}{a_\beta-a_\alpha}\right)
^{-\frac{1}{2}((D_\beta-D_\alpha)^2+c_1(X)(D_\beta-D_\alpha) )}\\
&&\cdot 
\prod_{\beta=1}^r g(a_\beta)^{-\frac{1}{2}
( D_\beta^2 + c_1(X)\cdot D_\beta)}\exp(-F_{X_0,A,B,\vD}(\va;\Lambda)).
\end{eqnarray*}
Therefore
$$
\log\left( \sum_{\sum D_\alpha =d}
 \Lambda^{|\vD|^2} h_\vD(\ep_1,\ep_2,\va;\Lambda)\right)
$$
is analytic in $\ep_1,\ep_2$ near $\ep_1=\ep_2=0$.
\end{proof}

By Lemma \ref{thm:logloglog}, the pole of
$\log Z^\inst_{X_0,A,B,d}$ along $\ep_1=\ep_2=0$
is the same as that of
\begin{eqnarray*}
&& -\log Z^\inst_{\bC^2,A,B}(w,u,\va;\Lambda)
 - \log Z^\inst_{\bC^2,A,B}(-w,u-kw,\va;\Lambda)\\
&=& \frac{\cF^\inst_{\bC^2,A,B}(w,u,\va;\Lambda) }{wu}
+ \frac{\cF^\inst_{\bC^2,A,B}(-w,u-kw,\va;\Lambda)}{-w(u-kw)}.
\end{eqnarray*}

\begin{df}[logarithm of the instanton part]
Define
$$
\cF^\inst_{X_0,d}(\ep_1,\ep_2,\va;\Lambda)
= -u(u-kw)\log Z^\inst_{X_0,d}(\ep_1,\ep_2,\va;\Lambda).
$$
\end{df}
\begin{thm} \label{thm:toric-Ctwo}
If $\cF^\inst_{\bC^2,A,B}(\ep_1,\ep_2,\va;\Lambda)$
is analytic in $\ep_1,\ep_2$ near $\ep_1=\ep_2=0$, then
\begin{enumerate}
\item[(a)] $\cF^\inst_{X_0,A,B,d}(\ep_1,\ep_2,\va;\Lambda)$ is analytic
in $\ep_1,\ep_2$ near $\ep_1=\ep_2=0$,
\item[(b)] $\displaystyle{\lim_{\ep_1,\ep_2\to 0} \cF^\inst_{X_0,A,B,d}(\ep_1,\ep_2,\va;\Lambda)
=k\lim_{\ep_1,\ep_2\to 0} \cF^\inst_{\bC^2,A,B}(\ep_1,\ep_2,\va;\Lambda)}$.
\end{enumerate}
\end{thm}
\begin{proof}
Let 
$$
g_k(w,u,\va;\Lambda) = -u(u-kw)\left(
\frac{\cF^\inst_{\bC^2,A,B}(w,u,\va;\Lambda) }{wu}
+\frac{\cF^\inst_{\bC^2,A,B}(-w,u-kw,\va;\Lambda)}{-w(u-kw)}
\right).
$$
Note that $(w,u)$ and $(\ep_1,\ep_2)$ are related by a coordinate
transformation in $SL(2,\bZ)$.
By Lemma \ref{thm:logloglog}, it suffices to show that
\begin{enumerate}
\item[(a)'] $g_k(w,u,\va;\Lambda)$ is analytic in $w$, $u$  near $w=u=0$,
\item[(b)'] $\displaystyle{\lim_{w,u\to 0}} g_k(w,u,\va;\Lambda)
=k\lim_{\ep_1,\ep_2\to 0} \cF^\inst_{\bC^2,A,B}(\ep_1,\ep_2,\va,\Lambda)$.
\end{enumerate}

We have
$$
\cF^\inst_{\bC^2}(-w,u-kw,\va;\Lambda)-\cF^\inst_{\bC^2,A,B}(w,u,\va;\Lambda) 
=w H_k(w,u,\va;\Lambda)
$$
where $H_k(w,u,\va;\Lambda)$ is analytic in $w,u$ near $w=u=0$.
So 
\begin{equation}\label{eqn:gFH}
g_k(w,u,\va;\Lambda) = k\cF_{\bC^2,A,B}^\inst(w,u,\va,\Lambda) + u H_k(w,u,\va;\Lambda).
\end{equation}

(a)' and (b)' are are immediate consequences of \eqref{eqn:gFH}.

\end{proof}
Theorem \ref{thm:Nekrasov-Ctwo} and Theorem \ref{thm:toric-Ctwo} imply:
\begin{thm}[Nekrasov  conjecture for toric surfaces: instanton part]
\label{thm:Nekrasov-toric}
$ $\\
\begin{enumerate}
\item {\em 4d pure gauge theory:}
\begin{enumerate}
\item $\cF^\inst_{X_0,d}(\ep_1,\ep_2,\va;\Lambda)$ is analytic
in $\ep_1,\ep_2$ near $\ep_1=\ep_2=0$.
\item $\displaystyle{\lim_{\ep_1,\ep_2\to 0} \cF^\inst_{X_0,d}(\ep_1,\ep_2,\va;\Lambda)}
= k\cF_0^\inst(\va,\Lambda)$, where
$\cF_0^\inst(\va,\Lambda)$ is the instanton part of the Seiberg-Witten prepotential
of 4d pure gauge theory.
\end{enumerate}

\item {\em 4d gauge theory with $N_f$ fundamental matter 
hypermultiplets:}
\begin{enumerate}
\item $\cF^\inst_{X_0,d}(\ep_1,\ep_2,\va,\vm;\Lambda)$ is analytic
in $\ep_1,\ep_2$ near $\ep_1=\ep_2=0$.
\item $\displaystyle{\lim_{\ep_1,\ep_2\to 0} \cF^\inst_{X_0,d}(\ep_1,\ep_2,\va,\vm;\Lambda)}
= k\cF_0^\inst(\va,\vm,\Lambda)$, where
$\cF_0^\inst(\va,\vm,\Lambda)$ is the instanton part of the Seiberg-Witten prepotential
of 4d gauge theory with $N_f$ fundamental matter hypermultiplets.
\end{enumerate}

\item {\em 4d gauge theory with one adjoint matter hypermultiplet:}
\begin{enumerate}
\item $\cF^\inst_{X_0,d}(\ep_1,\ep_2,\va,m;\Lambda)$ is analytic
in $\ep_1,\ep_2$ near $\ep_1=\ep_2=0$.
\item $\displaystyle{\lim_{\ep_1,\ep_2\to 0} \cF^\inst_{X_0,d}(\ep_1,\ep_2,\va,m;\Lambda)}
= k\cF_0^\inst(\va,m,\Lambda)$, where
$\cF_0^\inst(\va,m,\Lambda)$ is the instanton part of the Seiberg-Witten prepotential
of 4d gauge theory with one adjoint matter hypermultiplet.
\end{enumerate}

\item {\em 5d gauge theory compactified on a circle of circumference  $\beta$:}
\begin{enumerate}
\item $\cF^\inst_{X_0,d}(\ep_1,\ep_2,\va;\Lambda,\beta)$ is analytic
in $\ep_1,\ep_2$ near $\ep_1=\ep_2=0$.
\item $\displaystyle{\lim_{\ep_1,\ep_2\to 0} \cF^\inst_{X_0,d}(\ep_1,\ep_2,\va;\Lambda,\beta)}
= k\cF_0^\inst(\va,\Lambda,\beta)$, where
$\cF_0^\inst(\va,\Lambda,\beta)$ is the instanton part of the Seiberg-Witten prepotential
of 5d gauge theory compactified on a circle of circumference $\beta$.
\end{enumerate}
\end{enumerate}
\end{thm}

\section{The Perturbative Part}
\label{sec:perturbative}

In this section we prove the perturbative parts of the conjecture, 
of which instanton counterparts were proved in Theorem \ref{thm:Nekrasov-toric}.
The perturbative part comes from the difference between framed instantons
on the compact toric surface $X$ and 
unframed instantons on the noncompact toric surface $X_0$, 
so we must consider the virtual tangent and natural bundles
of the moduli space of unframed instantons on $X_0$. Evaluating 
the required multiplicative classes at such  bundles gives rise to 
infinite products which need to be regularised. Following \cite{NO} 
we use zeta-function regularization (Definition \ref{zeta}).

\subsection{The virtual tangent bundle of $\Mx$} \label{sec:Tpert}
Given $(E,\Phi)\in \MX$, we may look at $E|_{X_0}$ as representing a point in
the moduli space $\Mx$ of unframed instantons 
on the noncompact surface $X_0$. We have
\begin{eqnarray*}
&& \ch_{\tT}T^\vir_{E|_{X_0}} \Mx
=-\ch_{\tT} \Ex^*(E|_{X_0},E|_{X_0})\\
&=& \sum_{\alpha,\beta} e^{a_\beta-a_\alpha }
\sum_{v\in V(\Gamma)}  e^{w_{D_\beta}^v-w_{D_\alpha}^v}\Bigl(
N_{Y_\alpha^v,Y_\beta^v}(w_1^v,w_2^v) - \frac{1}{(1-e^{-w_1^v})(1-e^{-w_2^v})}\Bigr)\\
&=& \sum_{v\in \Gamma}
\sum_{\alpha,\beta} e^{ (a_\beta + w_{D_\beta}^v)- (a_\alpha+w_{D_\alpha}^v)}
\Bigl(N_{Y_\alpha^v,Y_\beta^v}(w_1^v,w_2^v) - \frac{1}{(1-e^{-w_1^v})(1-e^{-w_2^v})}\Bigr).
\end{eqnarray*}
The perturbative part of the $\tT$-equivariant Chern character 
of the tangent bundle is given by 
\begin{eqnarray*}
&& \ch_\tT T_{E|_{X_0}}^\pert \define \ch_{\tT} T^\vir_{E|_{X_0}}\Mx- \ch_{\tT}T_{(E,\Phi)}\MX\\
&=& -\sum_{\alpha,\beta} e^{a_\beta-a_\alpha} \Bigl(\frac{1}{(1-e^{-w})(1-e^u)}
+\frac{1}{(1-e^w)(1-e^{u-kw}) }\Bigr)\\
&=& \frac{-\sum_{\alpha,\beta} e^{a_\beta-a_\alpha}}{(1-e^u)(1-e^{u-kw}) }
\Bigl(1+ \sum_{j=1}^{k-1} e^{u-jw}\Bigr).
\end{eqnarray*}
\begin{ex} $X=\bP^2$, $X_0=\bC^2$.
\begin{eqnarray*}
\ch_\tT T_{(E,\Phi)}^\pert 
&=& -\sum_{\alpha,\beta} e^{a_\beta -a_\alpha} \Bigl(\frac{ 1}{(1-e^{\ep_2-\ep_1})(1-e^{-\ep_2})}
+\frac{1 }{(1-e^{\ep_1-\ep_2})(1-e^{-\ep_1}) }\Bigr) \\
&=&\frac{-\sum_{\alpha,\beta}  e^{a_\beta -a_\alpha}  }{ (1-e^{-\ep_1})(1-e^{-\ep_2}) }.
\end{eqnarray*}
\end{ex}

Let $A$ be a multiplicative class defined by a formal power series $f(x)$.
Formally, 
evaluating $A$ on the  
tangent bundle produces the following perturbative part:
\begin{equation}\label{eqn:ff}
A_{\tT}(T_{(E,\Phi)}^\pert)\  =
\frac{1}{\prod_{i,j=0}^\infty
f(a_\beta-a_\alpha-iw + ju) 
\prod_{i,j=0}^\infty
f(a_\beta-a_\alpha+ iw + j(u-kw)) }.
\end{equation}
The  infinite product on the right hand side  requires regularization.

\subsection{The natural virtual bundle}

Given $(E,\Phi)\in \MX$,  once again looking at 
$E|_{X_0}$ as representing a point in
$\Mx$, we have
\begin{eqnarray*}
&& \ch_{\tT}V^\vir_{E|_{X_0}}
=-\chi_{\tT} \Ex^* E\\
&=& \sum_{\beta} e^{a_\beta}
\sum_{v\in V(\Gamma)}  e^{w_{D_\beta}^v}\Bigl(
N_{Y_\beta^v}(w_1^v,w_2^v) - \frac{1}{(1-e^{-w_1^v})(1-e^{-w_2^v})}\Bigr)\\
&=& \sum_{v\in \Gamma}
\sum_{\alpha,\beta} e^{ (a_\beta + w_{D_\beta}^v)}
\Bigl(N_{Y_\beta^v}(w_1^v,w_2^v) - \frac{1}{(1-e^{-w_1^v})(1-e^{-w_2^v})}\Bigr).
\end{eqnarray*}
The perturbative part of the $\tT$-equivariant Chern character 
of the natural bundle is given by 
\begin{eqnarray*}
&& \ch_\tT V_{E|_{X_0}}^\pert \define \ch_{\tT} V^\vir_{E|_{X_0}}-  \ch_{\tT}V_{(E,\Phi)}\\
&=& -\sum_{\alpha,\beta} e^{a_\beta} \Bigl(\frac{1}{(1-e^{-w})(1-e^u)}
+\frac{1}{(1-e^w)(1-e^{u-kw}) }\Bigr)\\
&=& \frac{-\sum_{\beta} e^{a_\beta}}{(1-e^u)(1-e^{u-kw}) }
\Bigl(1+ \sum_{j=1}^{k-1} e^{u-jw}\Bigr).
\end{eqnarray*}
\begin{ex} $X=\bP^2$, $X_0=\bC^2$.
\begin{eqnarray*}
\ch_\tT V_{E|_{X_0}}^\pert 
&=& -\sum_{\beta} e^{a_\beta} \Bigl(\frac{ 1}{(1-e^{\ep_2-\ep_1})(1-e^{-\ep_2})}
+\frac{1 }{(1-e^{\ep_1-\ep_2})(1-e^{-\ep_1}) }\Bigr) \\
&=&\frac{-\sum_{\beta}  e^{a_\beta}  }{ (1-e^{-\ep_1})(1-e^{-\ep_2}) }.
\end{eqnarray*}
\end{ex}

Let $B$ be a multiplicative class defined by a formal power series $g(x)$.
Formally, evaluating $B$ on the natural  bundle produces the
following perturbative part:
\begin{equation}\label{eqn:gg}
B_{\tT}(V_{E|_{X_0}}^\pert)\  =
\frac{1}{\prod_{i,j=0}^\infty
g(a_\beta-iw + ju) 
\prod_{i,j=0}^\infty
g(a_\beta+ iw + j(u-kw)) }.
\end{equation}
The 
infinite product
on the right hand side
requires regularization.

\subsection{Regularization}

Following \cite[Appendix A]{NO}, we introduce the
following functions.

\begin{df}[zeta-regularization] \label{zeta}
\begin{equation}\label{eqn:g}
\gamma_{\ep_1,\ep_2}(x;\Lambda) \define \frac{d}{ds}
\Bigr|_{s=0}\frac{\Lambda}{\Gamma(s)} \int_0^\infty \frac{dt}{t} t^s
\frac{e^{-tx}}{(e^{\ep_1 t}-1)(e^{\ep_2 t}-1)}.
\end{equation}
\begin{equation}
\begin{aligned}
\gamma_{\ep_1,\ep_2}(x\mid \beta;\Lambda)\define&
\frac{1}{2\ep_1\ep_2}
\left(-\frac{\beta}{6} \bigl(x+\frac{1}{2}(\ep_1+\ep_2)\bigr)^3
+x^2 \log(\beta\Lambda)\right)\\
& +\sum_{n=1}^\infty \frac{1}{n}
\frac{e^{-\beta nx} }{(e^{\beta n\ep_1} -1)(e^{\beta n \ep_2}-1)}.
\end{aligned}
\end{equation}
\end{df}

$\exp(\gamma_{\ep_1,\ep_2}(x;\Lambda))$ is a regularization
of the infinite product
$$
\prod_{i,j=0}^\infty \frac{\Lambda}{x-i\ep_1-j\ep_2}.
$$

For a very nice explanation of this regularization scheme see 
\cite{O}.
The function $\gamma_{\ep_1,\ep_2}(x;\Lambda)$
satisfy the following properties (see \cite[Appendix A]{NO}):
\begin{fact}\label{gamma}
\begin{enumerate}
\item $\ep_1\ep_2\gamma_{\ep_1,\ep_2}(x;\Lambda)$
is analytic in $\ep_1$, $\ep_2$ near $\ep_1=\ep_2=0$;
\item $\displaystyle{\lim_{\ep_1,\ep_2\to 0} 
\ep_1\ep_2 \gamma_{\ep_1,\ep_2}(x;\Lambda) }
=-\frac{1}{2}x^2 \log\frac{x}{\Lambda} +\frac{3}{4} x^2$.
\end{enumerate}
\end{fact}

\subsection{Nekrasov conjecture:  perturbative part}

Applying zeta-regularization to \eqref{eqn:ff}
and \eqref{eqn:gg}, we obtain the following
definitions:

\begin{df}[perturbative part of the partition function]\label{df:Fpert}
$ $\\
\begin{enumerate}
\item 4d pure gauge theory:
\begin{eqnarray*}
&&\cF^\pert_{X_0, A=1, B=1}(\ep_1,\ep_2,\va;\Lambda)\\
&\define& u(u-kw)\cdot \biggl(\sum_{\alpha,\beta} (\gamma_{-w,u}(a_\beta-a_\alpha;\Lambda)
+ \gamma_{w,u-kw}(a_\beta-a_\alpha;\Lambda) )\biggr)
\end{eqnarray*}
$$
Z^\pert_{X_0, A=1, B=1}(\ep_1,\ep_2,\va;\Lambda)
\define \exp\left( \frac{\cF^\pert_{X_0, A=1, B=1}(\ep_1,\ep_2,\va;\Lambda)}{-u(u-kw)}\right),
$$
\item 4d gauge theory with $N_f$ fundamental matter hypermultiplets:
\begin{eqnarray*}
&& \cF^\pert_{X_0, A=1, B=E_{\vm}}(\ep_1,\ep_2,\va;\Lambda)\\
&\define & u(u-kw)\cdot \biggl(\sum_{\alpha,\beta} \bigl(\gamma_{-w,u}(a_\beta-a_\alpha;\Lambda)
+ \gamma_{w,u-kw}(a_\beta-a_\alpha;\Lambda\bigr) \\
&& -\sum_{\beta,f} \bigl(\gamma_{-w,u}(a_\beta +m_f;\Lambda)+
\gamma_{w,u-kw}(a_\beta+m_f,\Lambda)\bigr)\biggr)
\end{eqnarray*}
$$
Z^\pert_{X_0, A=1, B=E_{\vm}}(\ep_1,\ep_2,\va;\Lambda)
\define \exp\left( \frac{\cF^\pert_{X_0, A=1, B=E_{\vm}}(\ep_1,\ep_2,\va;\Lambda)}{-u(u-kw)}\right),
$$
\item 4d gauge theory with one adjoint matter hypermultiplet:\\
\begin{eqnarray*}
&& \cF^\pert_{X_0, A=E_m, B=1}(\ep_1,\ep_2,\va;\Lambda)\\
&\define & u(u-kw)\cdot \biggl(\sum_{\alpha,\beta} \bigl(\gamma_{-w,u}(a_\beta-a_\alpha;\Lambda)
-\gamma_{-w,u}(m+a_\beta-a_\alpha;\Lambda)\\
&& + \gamma_{w,u-kw}(a_\beta-a_\alpha;\Lambda) 
-\gamma_{w,u-kw}(m+a_\beta-a_\alpha;\Lambda\bigr)\biggr)
\end{eqnarray*}
$$
Z^\pert_{X_0, A=E_m, B=1}(\ep_1,\ep_2,\va;\Lambda)
\define \exp\left( \frac{\cF^\pert_{X_0, A=E_m, B=1}(\ep_1,\ep_2,\va;\Lambda)}{-u(u-kw)}\right),
$$
\item 5d gauge theory compactified at a circle of circumference $\beta$:
\begin{eqnarray*}
&& \cF^\pert_{X_0, A=\widehat{A}_\beta, B=1}(\ep_1,\ep_2,\va;\Lambda)\\
&\define& u(u-kw)\sum_{p,q} (\gamma_{-w,u}(a_p-a_q;\beta,\Lambda)
+ \gamma_{w,u-kw}(a_p-a_q;\beta,\Lambda) 
\end{eqnarray*}
$$
Z^\pert_{X_0, A=\widehat{A}_\beta, B=1}(\ep_1,\ep_2,\va;\Lambda)
\define \exp\left( \frac{\cF^\pert_{X_0, A=\widehat{A}_\beta, B=1}(\ep_1,\ep_2,\va;\Lambda)}{-u(u-kw)}\right)
.$$
\end{enumerate}
\end{df}

\begin{ex} $X=\bP^2$, $X_0=\bC^2$.
\begin{enumerate}
\item 4d pure gauge theory:
$$
\cF^\pert_{\bC^2, A=1, B=1}(\ep_1,\ep_2,\va;\Lambda)
= \ep_1 \ep_2 \sum_{\alpha,\beta} \gamma_{\ep_1,\ep_2}(a_\beta-a_\alpha;\Lambda),
$$
\item 4d gauge theory with $N_f$ fundamental matter hypermultiplets:
\begin{eqnarray*}
&& \cF^\pert_{\bC^2, A=1, B=E_{\vm}}(\ep_1,\ep_2,\va;\Lambda)\\
&= & \ep_1\ep_2 \biggl(
\sum_{\alpha,\beta} \gamma_{\ep_1,\ep_2}(a_\beta-a_\alpha;\Lambda)
-\sum_{\beta,f} \gamma_{\ep_1,\ep_2}(a_\beta +m_f;\Lambda)
\biggr),
 \end{eqnarray*}
\item 4d gauge theory with one adjoint matter hypermultiplet:
\begin{eqnarray*}
&& \cF^\pert_{\bC^2, A=E_m, B=1}(\ep_1,\ep_2,\va;\Lambda)\\
&= & \ep_1\ep_2 \sum_{\alpha,\beta}\biggl(
\gamma_{\ep_1,\ep_2}(a_\beta-a_\alpha;\Lambda)
-\gamma_{\ep_1,\ep_2}(m+a_\beta-a_\alpha;\Lambda\bigr)
\biggr),
\end{eqnarray*}
\item 5d gauge theory compactified at a circle of circumference $\beta$:
$$
\cF^\pert_{\bC^2, A=\widehat{A}_\beta, B=1}(\ep_1,\ep_2,\va;\Lambda)
= \ep_1 \ep_2 \sum_{p,q} \gamma_{\ep_1,\ep_2}(a_p-a_q\mid \beta;\Lambda).
$$
\end{enumerate}
\end{ex}

\begin{thm}[Nekrasov conjecture: perturbative part]\label{thm:toric-pert}
$ $\\
\begin{enumerate}
\item 4d pure gauge theory:
$$
\lim_{\ep_1,\ep_2\to 0} \cF^\pert_{X_0, A=1, B=1}(\ep_1,\ep_2,\va;\Lambda)\\
= k\cF^\pert_0(\va,\Lambda)
$$
where
$$
\cF^\pert_0(\va,\Lambda)= \sum_{\alpha\neq \beta}\left(
-\frac{1}{2}(a_\alpha-a_\beta)^2\log\left(\frac{a_\alpha-a_\beta}{\Lambda}\right) 
+\frac{3}{4}(a_\alpha-a_\beta)^2\right)
$$
is the perturbative part of the Seiberg-Witten prepotential of 4d pure gauge theory.

\item 4d gauge theory with $N_f$ fundamental matter hypermultiplets:
$$
\lim_{\ep_1,\ep_2\to 0}\cF^\pert_{X_0, A=1, B=E_{\vm}}(\ep_1,\ep_2,\va;\Lambda)
=k \cF_0^\pert(\va,\vm,\Lambda)
$$
where
\begin{eqnarray*}
\cF_0^\pert(\va,\vm,\Lambda)&=& \sum_{\alpha\neq \beta}
\Bigl(-\frac{1}{2}(a_\alpha-a_\beta)^2\log\left(\frac{a_\alpha-a_\beta}{\Lambda}\right) 
+\frac{3}{4}(a_\alpha-a_\beta)^2\Bigr)\\
&&  + \sum_{\beta,f}
\Bigl(\frac{1}{2}(a_\beta+ m_f)^2\log\left(\frac{a_\beta+m_f}{\Lambda}\right) 
-\frac{3}{4}(a_\beta+ m_f)^2\Bigr)
\end{eqnarray*}
is the perturbative part of the Seiberg-Witten prepotential of 4d gauge theory with
$N_f$ fundamental matter hypermultiplets.

\item 4d gauge theory with one adjoint matter hypermultiplet:\\
$$
\lim_{\ep_1,\ep_2\to 0}\cF^\pert_{X_0, A=E_m, B=1}(\ep_1,\ep_2,\va;\Lambda)
= k\cF_0^\pert(\va,m,\Lambda)
$$
where
\begin{eqnarray*}
\cF_0^\pert(\va,m,\Lambda)
&=&\sum_{\alpha\neq \beta}
\Bigl(-\frac{1}{2}(a_\alpha-a_\beta)^2\log\left(\frac{a_\alpha-a_\beta}{\Lambda}\right) 
+\frac{3}{4}(a_\alpha-a_\beta)^2\\
&&  +\frac{1}{2}(a_\alpha-a_\beta + m)^2\log\left(\frac{a_\alpha-a_\beta+m}{\Lambda}\right) 
-\frac{3}{4}(a_\alpha-a_\beta+ m)^2\Bigr)
\biggr)
\end{eqnarray*}
is the perturbative part of the Seiberg-Witten prepotential of 4d gauge
theory with one adjoint matter hypermultiplets.

\item 5d gauge theory compactified at a circle of circumference $\beta$.
$$
\lim_{\ep_1,\ep_2\to 0} \cF^\pert_{X_0, A=\widehat{A}_\beta, B=1}(\ep_1,\ep_2,\va;\Lambda)
= k \cF^\pert_0(\va,\Lambda,\beta)
$$
where
$$
\cF^\pert_0(\va,\Lambda,\beta)=\sum_{p\neq q}\left(-\frac{\beta}{12}(a_p-a_q)^3
+\frac{1}{2}(a_p-a_q)^2\log(\beta\Lambda)\right)
$$
is the perturbative part of the Seiberg-Witten prepotential of 5d gauge
theory compactified on a circle.
\end{enumerate}
\end{thm}

\begin{proof}
We prove (1), (2), (3). The proof of (4) is similar. 

Define
$$
f_k(u,w,x;\Lambda) = u(u-kw)(\gamma_{-w,u}(x;\Lambda) + \gamma_{w,u-kw}(x;\Lambda)).
$$  
By Definition \ref{df:Fpert} (definition of $\cF^\pert$), it suffices to show that
$$
\lim_{u,w\to 0} f_k(u,w,x;\Lambda) = 
k \left( -\frac{1}{2} x^2 \log\frac{x}{\Lambda} +\frac{3}{4} x^2\right).
$$
Let $g(\ep_1,\ep_2,x;\Lambda)= \ep_1\ep_2\gamma_{\ep_1,\ep_2}(x;\Lambda)$. Then by 
Fact \ref{gamma},
\begin{enumerate}
\item[(i)] $g(\ep_1,\ep_2,x;\Lambda)$ is analytic in $\ep_1,\ep_2$ near $\ep_1=\ep_2=0$.
\item[(ii)] $\displaystyle{\lim_{\ep_1,\ep_2\to 0} g(\ep_1,\ep_2,x;\Lambda) 
=-\frac{1}{2} x^2 \log\frac{x}{\Lambda} +\frac{3}{4}x^2}$.
\end{enumerate}
By (i), we have
$$
g(w,u-kw,x;\Lambda)- g(-w,u,x;\Lambda) = w h_k(u,w,x;\Lambda)
$$
where $h_k(u,w,x;\Lambda)$ is analytic in $w,u$ near $w=u=0$. We have
\begin{eqnarray*}
f_k(u,w,x;\Lambda) &=&  u(u-kw)\left(\frac{ g(-w,u,x;\Lambda) }{-w u } + \frac{g(w,u-kw;\Lambda)}{w(u-kw)}\right) \\
&=& k g(-w,u,x;\Lambda) + u h_k(u,w,x;\Lambda).
\end{eqnarray*}
Therefore
$$
\lim_{u,w\to 0} f_k(u,w,x;\Lambda) =  k\lim_{\ep_1,\ep_2\to 0} g(\ep_1,\ep_2,x;\Lambda)
= k \left( -\frac{1}{2} x^2 \log\frac{x}{\Lambda} +\frac{3}{4} x^2\right).
$$

\end{proof}

\begin{appendix}

\section{Kobayashi--Hitchin correspondence and 
existence of instantons}

 In this section we recall some  results relating instantons in 
pure gauge theory to
holomorphic bundles. 
 The  Kobayashi--Hitchin correspondence predicts an equivalence  between 
instantons and holomorphic bundles in various settings, see \cite{LT}.
 For an $SU(n)$  bundle $E$  over compact K\"ahler surface $X$ this 
correspondence was proved by Donaldson \cite{Do1}:
 The moduli space of irreducible anti-self-dual 
connections on $E$ is naturally identified with the set of 
equivalence classes of stable holomorphic $SL(n,\mathbb C)$ 
bundles which are topologically equivalent to $E$ (see \cite{DK} 
Corollary 6.1.6 for a proof of the rank 2 case). Note that here 
stability is taken with respect to the K\"ahler class. Under this 
correspondence the topological charge of the instanton 
corresponds to the second Chern number of the bundle.

 To obtain a Kobayashi--Hitchin correspondence
 over a  non-compact K\"ahler manifold $(X,\omega)$
one must impose some conditions on the 
behaviour of holomorphic bundles at infinity. 
The instanton charge is obtained by integration 
of the curvature of the connection over $X$, 
and the mildest constraint that guarantees 
finiteness of this integral is to demand 
that the curvature decays as $1/r^2.$ 

For a manifold $X$ that can be compactified to 
$\bar X = X \cup D$ by adding a
smooth divisor $D$ with positive normal bundle, 
Bando \cite{Ba} defined a notion on $U(r)$ flatness and 
proved the following: There is a correspondence  
between the moduli space of Hermitian--Einstein holomorphic vector 
bundles on $(X,\omega)$ whose curvature decays faster than 
$1/r^2$ with trivial holonomy at infinity and 
the moduli space of holomorphic vector bundles  $\bar X$
whose restriction to $D$ are $U(r)-$flat. 

Alternatively, one can study
non-compact Kobayashi-Hitchin correspondence
between  instantons and framed bundles, that is, 
holomorphic bundles that are trivialized at infinity. See 
 Donaldson \cite{Do2}  for first non-compact 
instance of the correspondence, namely  instantons on $\mathbb C^2$;
then King \cite{Ki} for instantons on the blow-up of 
$\mathbb C^2$; and Gasparim--K\"oeppe--Majumdar \cite{GKM} for 
instantons on $Z_k :=\mbox{Tot}{\mathcal O}_{\mathbb P^1}(-k).$

We remark that these correspondences refer to 
 classical  instantons, and corresponding non-compactified 
moduli spaces of  holomorphic 
vector bundles having $c_1=0$ (i.e. locally trivial sheaves), 
whereas in the supersymmetric case the vocabulary 
instanton moduli refers to the much more general notion of  
moduli of torsion free sheaves and their compactifications. 
In particular, existence of instantons with a prescribed  charge in
supersymmetric gauge theories can be obtained simply by 
considering non-locally free sheaves. Thus, existence 
results for supersymmetric instantons contrast
with existence of classical instantons, c.f.
 \cite{GKM}  Theorem 6.8, which says that the minimal local charge of a 
nontrivial $SU(2)$-instanton on $Z_k$ is $k-1.$

\section{Equivariant Cohomology}\label{sec:equivariant}
Let $ET$ be a contractible space on which
$T=(\bC^*)^k$ acts freely, and let $BT=ET/T$.
(For example, $ET=(\bC^\infty-\{0\})^k$
and $BT=(\bP^\infty)^k$.) Then
$ET\to BT$ is a universal principal 
$T$-bundle.

Suppose that $T=(\bC^*)^k$ acts on
an $m$-dimensional  complex manifold $M$. 
The $T$-equivariant
cohomology of $M$ is defined to be 
$$
H_T^*(M;\bQ) \define H^*(M_T;\bQ)
$$
where $M_T = M\times_T ET$.
There is a fibration $M_T \to BT=ET/T$
with fiber $M$. Let $i_M: M\to M_T$
be the inclusion of fiber. This induces a 
ring homomorphism
$$
i_M^*: H^*_T(M;\bQ)\to H^*(M;\bQ).
$$
In particular, when $M$ is a point,  the map
$$
i_\pt^*: H^*_T(\pt;\bQ)\cong \bQ[u_1,\ldots,u_k]
\to H^*(\pt;\bQ) \cong \bQ
$$
is given by $p(u_1,\ldots,u_k)\mapsto p(0,...,0)$,
where $u_1,\ldots,u_k \in H^2_T(\pt;\bQ)$.

\subsection{Integral}
Now suppose that $M$ is compact. 
Then integration along the fiber
gives $\bQ$-linear maps
\begin{equation}\label{eqn:intM}
\int_M: H^*(M;\bQ)\to H^*(\pt;\bQ)
\end{equation}
\begin{equation}\label{eqn:intMT}
\int_M: H^*_T(M;\bQ)=H^*(M_T;\bQ) \to H^*_T(\pt;\bQ)= H^*(BT;\bQ)
\end{equation}
such that
\begin{enumerate}
\item[(i)] $\int_M \alpha =0$ if $\alpha \in H^q(M;\bQ)$, $q<2m$.
\item[(ii)] $\int_M \alpha \in H^0(\pt)\cong\bQ$ if $\alpha\in H^{2m}(M;\bQ)$.
\item[(iii)] $\int_M\alpha =0$ if $\alpha \in H^q_T(M;\bQ)$, $q<2m$.
\item[(iv)] $\int_M\alpha \in H^{q-2m}_T(\pt;\bQ)$ if $\alpha \in H^q_T(M;\bQ)$, $q\geq 2m$.
Note that $H^{q-2m}_T(\pt;\bQ)=0$ when $q$ is odd, and $H^{q-2m}_T(\pt;\bQ)$ consists
of homogeneous polynomials in $u_1,\ldots,u_k$ of degree $q/2-m$ when $q$
is even.
\item[(v)] $i_\pt^* \int_M \alpha = \int_M i_M^* \alpha \in H^0(\pt;\bQ)\cong \bQ$
for $\alpha\in H_T^*(M;\bQ)$. 
\end{enumerate}

\subsection{Localization}
Let $M^T$ denote
the set of $T$-fixed points in $M$. Suppose
that each connected component 
of $M^T$ is a compact complex 
submanifold of $M$,  so that $M^T$ 
has a normal bundle
$N$ which is a complex vector
bundle. Note that $N$ might have
different ranks on different connected components of $M^T$.
$T$ acts on $M^T$ trivially, so $(M^T)_T = M^T\times BT$ and
$$
H_T^*(M^T;\bQ)\cong H^*(M^T;\bQ)\otimes_\bQ H_T(\pt;\bQ).
$$
The $T$-equivariant Euler class $e_T(N) \in 
H_T^*(M^T;\bQ)$ is invertible in
$$
H^*(M^T;\bQ)\otimes_\bQ \bQ[u_1,\ldots,u_k]_\m
$$
where $\bQ[u_1,\ldots,u_k]_\m$ is the localization
of the ring $\bQ[u_1,\ldots,u_k]$ at the maximal
ideal $\m$ generated by $u_1,...,u_k$.
The Atiyah-Bott localization formula says
\begin{equation}\label{eqn:localization}
\int_M \alpha =\int_{M^T}\frac{i^*\alpha}{e_T(N)}
\end{equation}
where $\alpha\in H^*_T(M;\bQ)$,
and $i^*: H_T^*(M;\bQ)\to H_T^*(M^T;\bQ)$
is induced by the inclusion $i: M^T\to M$.
In particular, if $M^T$ consists of isolated
points $p_1,\ldots,p_N$, then
\begin{equation}\label{eqn:points}
\int_M\alpha = \sum_{j=1}^N\frac{i_{p_j}^* \alpha}{e_T(T_{p_i}M)}
\end{equation}
where $i_{p_j}^*: H^*_T(M;\bQ)\to H^*_T(p_j;\bQ)\cong \bQ[u_1,\ldots,u_k]$
is induced by the inclusion $i_{p_j}: p_j \to M$.

Now suppose that $M$ is non-compact. Then
\eqref{eqn:intM} and \eqref{eqn:intMT} are  not defined.
However, when $M^T$ is  compact, we may {\em define}
\eqref{eqn:intMT} by the right hand side of \eqref{eqn:localization}.
Now (i), (ii), (v) are irrelevant, and (iii), (iv) do not hold:
given $\alpha\in H^q_T(M;\bQ)$, we have
$\int_M \alpha=0$ if $q$ is odd, and 
$\int_M\alpha$ is a rational function in $u_1,\ldots,u_k$
homogenous of degree $q/2-m$ (the degree can be 
negative).

\begin{ex}
Let $T_t=(\bC^*)^2$ act on $\bP^2$ by 
$(t_1,t_2)\cdot  [Z_0,Z_1,Z_2]= [Z_0, t_1 Z_1,t_2 Z_2]$.
We have $H^*_{T_t}(\pt;\bQ)= \bQ[\ep_1,\ep_2]$.
\begin{eqnarray*}
&& \int_{\bP^2} 1 =\frac{1}{\ep_1\ep_2} +\frac{1}{(-\ep_1)(-\ep_1+\ep_2)}
+\frac{1}{(-\ep_2)(\ep_1-\ep_2)} =0\\
&& \int_{\bC^2} 1 =\frac{1}{\ep_1\ep_2}
\end{eqnarray*}
\end{ex}

\subsection{Characteristic classes}
Let $c$ be a characteristic class for complex
vector bundles. Given a $T$-equivariant
complex vector bundle $V$ over $M$, 
$V_T = V\times_T ET$ is a vector bundle
over $M_T=M\times_T ET$. The $T$-equivariant
characteristic class $c_T$ is defined by
$$
c_T(E)\define c(E_T)\in H^*(M_T;\bQ)=H^*_T(M;\bQ).
$$

\section{Seiberg-Witten Prepotential} \label{sec:SW}

We  present a brief description of the Seiberg--Witten prepotential, 
which is described in detail  in the seminal work \cite{SW}, where Seiberg and
Witten gave an exact solution to $N=2$ supersymmetric Yang--Mills
in 4 dimensions with group $SU(2).$ For more details see also
\cite{NY} and \cite{Do}. For gauge theory with matter see 
\cite{DW} and \cite{BFMT}. The subject of 5d gauge theories
compactified on a circle and the corresponding Seiberg-Witten curves 
were introduced in \cite{Ne1}.

\subsection{$SU(2)$ case}   
The constraints of $N=2$ SUSY imply
that the quantum moduli space is the same as the classical one as
an algebraic variety.
Basic quantities are then the  coordinates
$u$ of the moduli space and the electric charge $a$, which in the
classical theory are related simply by $u=a^2/2$; in the quantum
theory this relation holds approximately for $u\rightarrow \infty$
by asymptotic freedom, but for finite $u$ the relation is much more
intricate and encodes fundamental geometric and physical information.
The description of the theory via the low energy effective Lagrangian
presents measurable quantities  as functions of the coordinates $u$
of the moduli space, and in
particular the electric charge $a = a(u).$ Moreover, Seiberg
\cite{Se} shows that the magic  of
supersymmetry allows the effective Lagrangian to
be expressed in terms of a single locally defined meromorphic function:
the prepotential $\cF_0$; all remaining quantities in the theory
being expressible as functions of $\cF_0$ and $a$.
An appropriate incarnation of Montonen--Olive duality accounts for
the appearance of the dual variable
$$a^D = \frac{d \cF_0}{da}$$
whose physical meaning is of the dual, that is, magnetic charge.
The defining relations giving
$$\tau= \frac{da^D}{da}, \,\,\, \tau^D = \frac{d(-a)}{da^D},$$
which imply that the duality transformation is
$$\tau^D = - \tau (a)^{-1}$$ and specializes to the Montonen--Olive
transformation $g^D= g^{-1}$ when the phase angle
$\theta = 0,$ but not otherwise.
The moduli space then  acquires expressions for  a K\"ahler  metric
$$ds^2 = Im (\tau da d\bar a)$$
with K\"ahler potential $\sum \frac{d\cF_0}{da_i} \bar a_i$, where $\tau$
is the matrix of periods $$\tau = \frac{d^2 \cF_0}{da^2}= \frac{da^D}{da}.$$
For $SU(2)$ the low-energy effective values of this coupling are given by
$\tau = \frac{\theta}{2\pi}+ \frac{4\pi i}{g^2}$ where $\theta$ is
is defined only modulo $2\pi {\mathbb Z}$;
consequently $\tau$ is defined only modulo $\mathbb Z$ and there is a
second
transformation fixing $a$ and taking $\tau \mapsto \tau + 1$.
Since $\tau = \frac{da^D}{da}$, it follows that $a^D\mapsto a^D+a.$
This pair of transformations acts as multiplication on the $2-$vector 
$(a^D,a)$
by the matrices
$$\left(\begin{matrix}0 & 1 \cr -1 & 0 \end{matrix}\right)\,\, \mbox{and}
\left(\begin{matrix}1 & 1 \cr 0 & 1 \end{matrix}\right)
$$  and fractional-linearly on $\tau$, 
thus  generating an $SL(2,\mathbb Z)$ action.
The upshot is that what lives intrinsically over a point $u$ in the
moduli space is not the electric charge $a(u)$ but the unimodular
lattice $\mathbb Z a(u)  + \mathbb Z a^D(u)$ of all electric and magnetic
charges. As $u$ varies we obtain a $\mathbb Z^2$ local system $V$ over
the moduli space, which Seiberg and Witten showed to have as simple as
possible behaviour; thus having
only 3 singularities at $\pm 1 $ and $\infty.$
Fixing a section of $V$
determines the prepotential up to a constant.
>From a careful analysis of  the monodromies at the singular points,
it follows that the local system itself can be identified
with the fiber cohomology of the elliptic curve
$$E_u\colon y^2 = (x+1)(x-1) (x-u).$$
The complexification $V_{\mathbb C}$ can be globally trivialized
in terms of a holomorphic 1-form $\lambda_1= \frac{dx}{y}$ and
a residueless meromorphic form $\lambda_2=\frac{xdx}{y}.$
One then chooses a homology basis consisting of a loop
$\gamma$ around the branch points $0,1$ and a loop $\gamma^D$ around
$1,u$; and  using such a basis, the correct geometric solution for the period is 
$$
\tau_u = \frac{ \oint_{\gamma^D}\lambda_1}{\oint_{\gamma} \lambda_1}.$$
In this solution,  $a$ and $a^D$ appear as the periods of   $\gamma$ and
$\gamma^D$  of the meromorphic 1-form
$$\lambda = \frac{ydx}{x^2-1}=\frac{(x-u)dx}{y}= \lambda_2-u\lambda_1.$$

\subsection{Higher rank case}

The Seiberg Witten solution is sometimes presented 
in reverse order, starting directly with the family of curves parametrized
by $u$ as we just described.
For instance, the solution for the group $SU(r)$ then appears as follows.
Let $\phi$ be an $SU(r)$ gauge field. Then
$$\det(xI-\phi) = x^r + U_2 x^{r-2} -U_3 x^{r-3} + ... +(-1)^r U_r,$$
where $U_k$ is the elementary symmetric polynomial of the eigenvalues
of $\phi$, with $U_1=0$ because $\phi$ takes values in $SU(r)$.
These are gauge invariant operators, so their
vacuum expectation values $u_k= \langle U_k \rangle$
 serve as coordinates of the classical moduli space.
These are the coordinates on the $\vu$-space: $u_2, ..., u_r$, 
which generalises the so-called $u$-plane in 
the $SU(2)$ case.

In case of added matter, then the duality transformations take a
different form, e.g. adding $N_f$  fundamental matter
hypermultiplets, the duality transformation becomes:
$$
{a^D \choose a} \mapsto R {a^D \choose a} + \sum_{i=1}^{N_f}m_i {n_i^D\choose n_i}
$$
where $R \in Sp(2(r-1), \mathbb Z),$ the $m_i$ are the masses of the $N_f$ 
particles added, and $n_i,n_i^D$ are integral $r\times r$ matrices.
Correspondingly, on the total space of the family of curves, 
there are then $N_f$ divisors ${\mathcal D}_i$ along which 
the meromorphic  differential $\lambda$ acquires a pole 
with constant residue $\frac{m_i}{2\pi \sqrt{-1}}$. Here again the charges
$a,a^D$ can be recovered as the periods of $\lambda$ over 
$\gamma$ and $\gamma^D.$

We now describe the Seiberg-Witten prepotential in various
gauge theories with gauge group $SU(r)$, starting directly with the 
Seiberg--Witten curves.
Consider the family of hyperelliptic curves of genus $r-1$ parametrized by
$\Lambda$, $\vu=(u_2,\ldots,u_r)$, and possibly some extra parameters,
in the following cases:
\begin{enumerate}
\item {\em 4d pure gauge theory} (see e.g. \cite[(4.5)]{NO}):
$$
C_\vu: \Lambda^{r}(w+\frac{1}{w})=P(z)= z^r + u_2 z^{r-2}+ \cdots + u_r .
$$

\item {\em 4d gauge theory with $N_f$ fundamental matter hypermultiplets}
(see e.g. \cite[(1.10)]{Ne2}):
$$
C_{\vu,\vm}: \Lambda+\frac{\Lambda^{2r-N_f}Q(z)}{w}=P(z),\quad
Q(z)=\prod_{f=1}^{N_f}(z-m_f).
$$

\item {\em 4d gauge theory with adjoint matter hypermultiplets} (see e.g. \cite[(6.32)]{NO}):
in this case the SW curve is the spectral curve of the elliptic 
Calogero--Moser system.
$$
C_{\vu, m}: \mathrm{Det}_{l,n}(L(w) - z)  = 0 , 
$$
where
$$
L_{l,n}(w) = \delta_{ln}\left(p_n+ \frac{m}{2\pi \sqrt{-1}}\log(\theta_{11}(w))'
\right)+ \frac{m}{2\pi \sqrt{-1}}(1-\delta_{ln})\frac{\theta_{11}(w+q_l-q_n)\theta_{11}'(0)}{\theta_{11}(w)\theta_{11}(q_l-q_n)}. 
$$
$$
\theta_{11}(\varpi;\tau)= \sum_{n\in \bZ}
e^{\pi\sqrt{-1} \tau(n+\frac{1}{2})^2 + 2\pi\sqrt{-1} (\varpi+\frac{1}{2})(n+\frac{1}{2}) }.
$$

\item {\em 5d gauge theory compactified at a circle of circumference $\beta$}
(see e.g. \cite[(7.19)]{NO}):
$$
C_{\vu,\beta}: (\beta\Lambda)^r (w+\frac{1}{w})
=X^{-r/2}P(X),\quad X= e^{\beta z}.
$$
\end{enumerate}

The {\em Seiberg-Witten differential} is
$$
dS=\frac{1}{2\pi\sqrt{-1}} z \frac{dw}{w}= 
\frac{1}{2\pi\sqrt{-1}}\frac{z P'(z) dz}{y}.
$$

Let $\{A_\alpha,\ B_\beta \mid \alpha,\beta=2,\ldots,r\}$
be a symplectic basis of $H_1(C_\vu,\bZ)$. Define functions
$a_\alpha$, $a^D_\beta$ on the $\vu$-plane by
$$
a_\alpha=\oint_{A_\alpha} dS,\quad
a^D_\alpha =2\pi\sqrt{-1} \oint_{B_\beta} dS.
$$
Then
$$
\omega_p= \frac{1}{2\pi\sqrt{-1}}\frac{z^{r-p}dz}{y},\quad
p=2,\ldots,r
$$
form a basis of holomorphic differentials on $C_\vu$.
The period matrix $\tau=(\tau_{\alpha\beta})$ is given by
$$
\tau_{\alpha\beta} =\frac{1}{2\pi{\sqrt{-1}}} \frac{\pa a^D_\alpha}{\pa a_\beta}.
$$
Note that a change of symplectic basis corresponds
to an element in $Sp(2(r-1),\bZ)$, the group of
duality  acting on the
period matrix $\tau=(\tau_{\alpha\beta}).$  In the $SU(2)$ or $U(2)$ cases,
we have  $r=2$, so  the group of duality is $Sp(2,\bZ)=SL(2,\bZ)$
and the SW curve is an elliptic curve.

The {\em Seiberg-Witten prepotential} is a locally defined function 
satisfying 
$$
a_{\alpha}^D = \frac {\partial \cF_0} {\partial a_{\alpha}}.
$$
Therefore the Seiberg-Witten prepotential and the peroid matrix
are related by
$$
\tau_{\alpha\beta} = \frac{1}{2\pi\sqrt{-1}} 
\frac{\pa^2 \cF_0}{\pa a_\alpha \pa a_\beta}.
$$

The full Seiberg--Witten prepotential is expressed as a sum
$$
\cF_0 = \cF_0^\pert + \cF_0^\inst
$$ 
where $\cF_0^\pert$ is the  {\em perturbative part}
and $\cF_0^\inst$ is the {\em instanton part}. 
The explicit expressions of the perturbative parts $\cF_0^\pert$ 
of the SW prepotentials in gauge theories (1), (2),
(3), (4) are given explicitly in (1), (2), (3), (4) of
Theorem \ref{thm:toric-pert}, respectively;
they have logrithm singularities along $\Lambda=0$. 
The instanton part $\cF_0^\inst$ of the SW prepotential
is a power series in $\Lambda^{2r}$: 
$$
\cF_0^\inst = O(\Lambda^{2r}) = f_1\Lambda^{2r} + f_2\Lambda^{4r} +
\cdots + f_n\Lambda^{2nr}+ \cdots 
$$
The coefficient  $f_n$ coming from the
$n$-instanton moduli space is
called the {\em $n$-th instanton correction} 
to the prepotential.

For further details we refer to \cite{DW}, \cite{GNY2}, \cite{Ne1},  
\cite{NO}, and \cite[Section 2]{NY}.

\end{appendix}

\end{document}